\documentclass[a4paper]{amsart}

\RequirePackage{amsmath}
\RequirePackage{bm}
\RequirePackage{amssymb}
\RequirePackage{upref}
\RequirePackage{amsthm}
\RequirePackage{enumerate}
\RequirePackage{pb-diagram}
\RequirePackage{amsfonts}
\RequirePackage[mathscr]{eucal}
\RequirePackage{verbatim}
\RequirePackage{xr}
\RequirePackage{graphicx}
\usepackage{calc}
\usepackage{xspace}
\RequirePackage{color}
\RequirePackage{ifthen}

\newcommand{\cf}{cf.\@\xspace}
\newcommand{\resp}{resp.\@\xspace}

\newcommand{\al}{\alpha}
\newcommand{\bet}{\beta}
\newcommand{\ga}{\gamma}
\newcommand{\de}{\delta }
\newcommand{\e}{\epsilon}

\newcommand{\f}{\varphi}
\newcommand{\h}{\eta}

\newcommand{\ka}{\kappa}
\newcommand{\lam}{\lambda}

\newcommand{\n}{\nu}

\newcommand{\vt}{\vartheta}

\newcommand{\s}{\sigma}
\newcommand{\x}{\xi}

\newcommand{\C}{\varGamma}

\newcommand{\F}{\varPhi}
\newcommand{\Lam}{\varLambda}

\newcommand{\di}[1]{#1\nobreakdash-\hspace{0pt}dimensional}

\newcommand{\fv}[2]{#1\hspace{0pt}_{|_{#2}}}
\newcommand{\cchi}[1]{\chi\hspace{0pt}_{_{#1}}}
\newcommand{\tchi}[1]{\tilde\chi\hspace{0pt}_{_{#1}}}
\newcommand{\so}{{\mc S_0}}

\newcommand{\const}{\tup{const}}

\newcommand{\msp[1]}[1]{\mspace{#1mu}}

\newcommand{\R}[1][n+1]{{\protect\mathbb R}^{#1}}

\newcommand{\Hh}[1][n+1]{{\protect\mathbb H}^{#1}}
\newcommand{\Ss}[1][n+1]{{\protect\mathbb S}^{#1}}

\newcommand{\N}{{\protect\mathbb N}}

\newcommand{\Z}{{\protect\mathbb Z}}
\newcommand{\eR}{\stackrel{\lower1ex \hbox{\rule{6.5pt}{0.5pt}}}{\msp[3]\R[]}}
\newcommand{\eN}{\stackrel{\lower1ex \hbox{\rule{6.5pt}{0.5pt}}}{\msp[1]\N}}
\newcommand{\eO}{\stackrel{\lower1ex \hbox{\rule{6pt}{0.5pt}}}{\msc O}}

\DeclareMathOperator{\graph}{graph}

\newcommand\im{\implies}
\newcommand\ra{\rightarrow}

\newcommand\pa{\partial}
\newcommand\pde[2]{\frac {\partial#1}{\partial#2}}
\newcommand\pd[3]{\frac {\partial#1}{\partial#2^#3}}   
\newcommand\pdc[3]{\frac {\partial#1}{\partial#2_#3}}   
\newcommand\pdm[4]{\frac {\partial#1}{\partial#2_#3^#4}}   
 
\newcommand{\un}{\infty}
\newcommand{\A}{\forall}

\newcommand{\set}[2]{\{\,#1\colon #2\,\}}
\newcommand{\uu}{\cup}
\newcommand{\ii}{\cap}
\newcommand{\uuu}{\bigcup}

\newcommand{\uud}{ \stackrel{\lower 1ex \hbox {.}}{\uu}}
\newcommand{\uuud}[1]{ \stackrel{\lower 1ex \hbox {.}}{\uuu_{#1}}}
\newcommand\su{\subset}

\newcommand\eS{\emptyset}
\newcommand{\sminus}[1][28]{\raise 0.#1ex\hbox{$\scriptstyle\setminus$}}

\newcommand{\wed}{\wedge}

\newcommand\ti{\times }

\newcommand{\abs}[1]{\lvert#1\rvert}

\newcommand{\norm}[1]{\lVert#1\rVert}

\newcommand{\spd}[2]{\protect\langle #1,#2\protect\rangle}

\newcommand\ch[3]{\varGamma_{#1#2}^#3}
\newcommand\cha[3]{{\bar\varGamma}_{#1#2}^#3}

\newcommand{\riem}[4]{R_{#1#2#3#4}}
\newcommand{\riema}[4]{{\bar R}_{#1#2#3#4}}

\newcommand{\tit}{\textit}

\newcommand{\tup}{\textup}

\newcommand{\mc}{\protect\mathcal}
\newcommand{\msc}{\protect\mathscr}

\providecommand{\bysame}{\makebox[3em]{\hrulefill}\thinspace}

\newcommand{\cq}[1]{\glqq{#1}\grqq\,}

\newcommand\fr[2]{\frac{#1}{1-\frac14{#2}^2}}

\newcommand{\bt}{\begin{thm}}
\newcommand{\bl}{\begin{lem}}
\newcommand{\bc}{\begin{cor}}
\newcommand{\bd}{\begin{definition}}
\newcommand{\bpp}{\begin{prop}}
\newcommand{\br}{\begin{rem}}
\newcommand{\bn}{\begin{note}}
\newcommand{\be}{\begin{ex}}
\newcommand{\bes}{\begin{exs}}
\newcommand{\bb}{\begin{example}}
\newcommand{\bbs}{\begin{examples}}
\newcommand{\ba}{\begin{axiom}}
\newcommand{\bas}{\begin{assumption}}

\newcommand{\et}{\end{thm}}
\newcommand{\el}{\end{lem}}
\newcommand{\ec}{\end{cor}}
\newcommand{\ed}{\end{definition}}
\newcommand{\epp}{\end{prop}}
\newcommand{\er}{\end{rem}}
\newcommand{\en}{\end{note}}
\newcommand{\ee}{\end{ex}}
\newcommand{\ees}{\end{exs}}
\newcommand{\eb}{\end{example}}
\newcommand{\ebs}{\end{examples}}
\newcommand{\ea}{\end{axiom}}
\newcommand{\eas}{\end{assumption}}

\newcommand{\bp}{\begin{proof}}
\newcommand{\ep}{\end{proof}}
\newcommand{\eps}{\renewcommand{\qed}{}\end{proof}}

\newcommand{\bal}{\begin{align}}

\newcommand{\bi}[1][1.]{\begin{enumerate}[\upshape #1]}
\newcommand{\bia}[1][(1)]{\begin{enumerate}[\upshape #1]}
\newcommand{\bin}[1][1]{\begin{enumerate}[\upshape\bfseries #1]}
\newcommand{\bir}[1][(i)]{\begin{enumerate}[\upshape #1]}
\newcommand{\bic}[1][(i)]{\begin{enumerate}[\upshape\hspace{2\cma}#1]}
\newcommand{\bis}[2][1.]{\begin{enumerate}[\upshape\hspace{#2\parindent}#1]}
\newcommand{\ei}{\end{enumerate}}

\newcommand\ndots{\raise 0.47ex \hbox {,}\hskip0.06em\cdots %
     \raise 0.47ex \hbox {,}\hskip0.06em} 

\newcommand{\q}{\quad}
\newcommand{\qq}{\qquad}

\newcommand{\hp}{\hphantom}

\newcommand\nd{\noindent}

\newskip\Csmallskipamount                                                
\Csmallskipamount=\smallskipamount
\newskip\Cmedskipamount
\Cmedskipamount=\medskipamount
\newskip\Cbigskipamount
\Cbigskipamount=\bigskipamount

\newcommand\cvs{\vspace\Csmallskipamount}   
\newcommand\cvm{\vspace\Cmedskipamount}

\newskip\csa
\csa=\smallskipamount

\newskip\cma
\cma=\medskipamount

\newskip\cba
\cba=\bigskipamount

\newdimen\spt
\newcommand\citem{\cvs\advance\itemno by
1{(\romannumeral\the\itemno})\hskip3pt}
\newcommand{\bitem}{\cvm\nd\advance\itemno by
1{\bf\the\itemno}\hspace{\cma}}

\newcount\itemno
\newcommand{\lae}[1]{\label{E:#1}}

\newcommand{\lal}[1]{\label{L:#1}}
\newcommand{\lad}[1]{\label{D:#1}}
\newcommand{\lac}[1]{\label{C:#1}}

\newcommand{\rl}[1]{Lemma~\ref{L:#1}}
\newcommand{\rd}[1]{Definition~\ref{D:#1}}

\newcommand{\rc}[1]{Corollary~\ref{C:#1}}

\newcommand{\re}[1]{\eqref{E:#1}}

\newcommand{\frl}[1]{Lemma~\ref{L:#1} on page~\tup{\pageref{L:#1}}}

\newcommand{\fre}[1]{\eqref{E:#1} on page~\tup{\pageref{E:#1}}}

\newskip\thmskip
\thmskip=\parindent

\newskip\hsk
\newenvironment{hinw}{\labelsep=0pt\begin{list}{}{\labelsep=0pt\itemindent=0pt\labelwidth=0pt\leftmargin=\parindent\rightmargin=0pt\partopsep=\cba}%
\item\it\nopagebreak\nopagebreak}%
{\end{list}}

\newcommand\bh{\begin{hinw}}
\newcommand{\eh}{\end{hinw}}

\newtheoremstyle{normal}
  {\cba}
  {\cba}
  {}
  {\thmskip}
  {\bfseries}
  {.}
  {\hsk}
  {}

\newtheoremstyle{abschnitt}
  {\cba}
  {\cba}
  {}
  {\thmskip}
  {\bfseries}
  {.}
  {\hsk}
  {}

\newtheoremstyle{italic}
  {\cba}
  {\cba}
  {\itshape}
  {\thmskip}
  {\bfseries}
  {.}
  {\hsk}
  {}

\newtheoremstyle{aufgaben}
  {\cba}
  {\cba}
  {}
  {}
  {\normalsize\bfseries}
  {.}
  {\hsk}
  {}

\newtheoremstyle{break}
  {\cba}
  {\cba}
  {\itshape}
  {}
  {\bfseries}
  {.}
  {\newline}
  {}

\swapnumbers
\theoremstyle{italic}
\newtheorem{thm}[subsection]{Theorem}
\newtheorem{lem}[subsection]{Lemma}
\newtheorem{prop}[subsection]{Proposition}
\newtheorem{cor}[subsection]{Corollary}

\theoremstyle{normal}
\newtheorem{rem}[subsection]{Remark}
\newtheorem{definition}[subsection]{Definition}
\newtheorem{example}[subsection]{Example}
\newtheorem{examples}[subsection]{Examples}
\newtheorem{ex}[subsection]{Exercise}
\newtheorem{note}[subsection]{}
\newtheorem{axiom}[subsection]{Axiom}
\newtheorem{assumption}[subsection]{Assumption}

\theoremstyle{aufgaben}
\newtheorem{exs}[subsection]{Exercises}

\swapnumbers

\numberwithin{equation}{section}
\numberwithin{figure}{section}

\newenvironment{textequation}[1][0.8]
{\begin{equation}
\begin{aligned}
\begin{minipage}{#1\linewidth}}
{\end{minipage}
\end{aligned}
\end{equation}
\ignorespacesafterend}

\newcommand{\btext}{\begin{textequation}}
\newcommand{\etext}{\end{textequation}}

\def\hinweis{\@startsection{subsection}{2}%
 \z@{0.7\linespacing\@plus 0.5\linespacing}{0.7\linespacing}%
{\normalfont\itshape\indent}}

\newcounter{hours}\newcounter{minutes}
\newcommand{\printtime}{%
\setcounter{hours}{\time/60}%
\setcounter{minutes}{\time-\value{hours}*60}%
\ifthenelse{\value{minutes}<9}{\thehours :0\theminutes}{\thehours:\theminutes}}

\usepackage[german,english]{babel}
\usepackage{graphicx}
\RequirePackage{amsmath}
\RequirePackage{bm}
\RequirePackage{amssymb}
\RequirePackage{upref}
\RequirePackage{amsthm}
\RequirePackage{enumerate}
\RequirePackage{pb-diagram}
\RequirePackage{amsfonts}
\RequirePackage[mathscr]{eucal}
\RequirePackage{verbatim}
\RequirePackage{xr}
\RequirePackage{graphicx}
\usepackage{calc}
\usepackage{xspace}

\makeatletter
\RequirePackage{color}
\newcommand{\ann}[1]{\renewcommand{\@makefnmark}{\mbox{$^{\color{red}{\@thefnmark}}$}}%
\footnote {#1}}
\makeatother

\RequirePackage{upref}
\RequirePackage{amsthm}
\RequirePackage{enumerate}
\usepackage[mathscr]{eucal}

\usepackage{xr-hyper}

\usepackage{calc}

\newlength{\oddsidemarginlength}
\newlength{\topmarginlength}

\newcounter{numberoflines}
\newcounter{tempcc}
\oddsidemargin=\oddsidemarginlength
\evensidemargin=\oddsidemargin
\topmargin=\topmarginlength
\usepackage[colorlinks=true,linkcolor=blue,citecolor=blue,urlcolor=blue]{hyperref}  

\begin{document}

\flushbottom


\title{Inverse Curvature Flows in Hyperbolic Space}

\author{Claus Gerhardt}
\address{Ruprecht-Karls-Universit\"at, Institut f\"ur Angewandte Mathematik,
Im Neuenheimer Feld 294, 69120 Heidelberg, Germany}
\email{gerhardt@math.uni-heidelberg.de}
\urladdr{\href{http://web.me.com/gerhardt/}{http://web.me.com/gerhardt/}}
\thanks{This work has been supported by the DFG}
\dedicatory{Dedicated to Stefan Hildebrandt on the occasion of his 75th birthday.}

%
\date{\today}
%


\begin{abstract} 
We consider inverse curvature flows in $\Hh$  with star-shaped initial hypersurfaces and prove that the flows exist for all time, and that the leaves converge to infinity, become strongly convex exponentially fast and also more and more totally umbilic.  After an appropriate rescaling the leaves converge in $C^\infty$ to a sphere.
\end{abstract}

\maketitle

\tableofcontents

\setcounter{section}{0}
\section{Introduction}
Curvature flows (driven by extrinsic curvatures) of compact hypersurfaces in a Riemannian space generally exist only for a finite time and then develop a singularity provided the flow is a pure curvature flow without an additional force term. This phenomenon occurs in the case of direct flows, which can also be characterized as \tit{contracting} flows, \cf \cite{gh:mean}, as well as for inverse flows, which can also be characterized as \tit{expanding} flows, see \cite{hi:penrose}.

In non-compact spaces of constant curvature we can expect that the inverse flows behave differently than the direct flows, since the inverse flows of geodesic spheres exist for all time. In \cite{cg90} we proved that inverse curvature flows of star-shaped hypersurfaces in Euclidean space exist for all time, converge to infinity and, after rescaling, converge to spheres.

In this paper we want to prove a similar result in hyperbolic space $\Hh$, $n\ge 2$. The initial hypersurface $M_0$ is supposed to be star-shaped with respect to a given point $p\in\Hh$, i.e., after introducing geodesic polar coordinates with center $p$, $M_0$ can be written as a graph over a geodesic sphere with center $p$ which we identify topologically with the standard sphere $\Ss[n]$. Let $F$ be a smooth curvature function,   homogeneous of degree $1$, monotone, and concave, defined in a symmetric, convex, open cone $\C\su\R[n]$, such that
\begin{equation}
\fv F\C>0\q\wed\q\fv F{\pa\C}=0.
\end{equation}
Then we consider the curvature flow
\begin{equation}
\dot x=-\F\nu
\end{equation}
with initial hypersurface $M_0$, the principal curvatures of which are supposed to lie in the cone $\C$; such a hypersurface is called \tit{admissible}. Here, the function $\F$ is defined by
\begin{equation}\lae{1.3}
\F=\F(r)=-r^{-1},\qq r>0,
\end{equation}
and $\F\nu$ stands for
\begin{equation}
\F(F)\nu,
\end{equation}
i.e., the flow equation is
\begin{equation}\lae{1.5}
\dot x=\frac 1F\nu,
\end{equation}
where $\nu$ is outward normal.

However, to simplify comparisons with former results and formulas, and also to make a generalization to more general flows easier, most results in this paper are formulated, and some are also proved, for a general smooth real valued function $\F$ defined on the positive real axis satisfying
\begin{equation}
\dot\F>0\q\wed\q \ddot\F\le 0.
\end{equation}

We shall normalize $F$ such that
\begin{equation}
F(1,\ldots,1)=n
\end{equation}
and shall also use the same notation $F$ when we assume $F$ to depend on the second fundamental form $h^i_j$ instead of the principal curvatures. Sometimes we also use the notation $\breve h^i_j$ for the second fundamental form of a hypersurface embedded in $\Hh$ to distinguish it from the second fundamental form $h^i_j$ of the same hypersurface viewed as being embedded in $\R$, which will happen when we parameterize $\Hh$ over the open ball $B_2(0)\su\R$.

We can now state our first result.
\bt
The flow \re{1.5} with a smooth and admissible initial hypersurface $M_0$ exists for all time. The flow hypersurfaces in hyperbolic space converge to infinity,  become strongly convex exponentially fast and also more and more totally umbilic. In fact there holds
\begin{equation}
\abs{\breve h^i_j-\de^i_j}\le c e^{-\frac tn},
\end{equation}
i.e., the principal curvatures are uniformly bounded and converge exponentially fast to $1$.
\et

For a more detailed analysis of the asymptotic behaviour we have to parameterize $\Hh$ over $B_2(0)\su \R$ such that the metric can be expressed in the form
\begin{equation}
\begin{aligned}
d\breve s^2&=\frac1{(1-\frac14\abs x^2)^2}dx^2\\
&=\frac1{(1-\frac14 r^2)^2}\{dr^2+r^2\s_{ij}dx^idx^j\}.
\end{aligned}
\end{equation}

The flow hypersurfaces $M(t)$ can now also be viewed as graphs
\begin{equation}
M(t)=\graph u(t,\cdot)
\end{equation}
over $\Ss[n]$ in Euclidean space, such that $0<u<2$, and convergence to infinity is tantamount to $u\ra 2$. The second fundamental form in $\R$ is denoted by $h^i_j$, or simply by $A$, where we omit the tensor indices. Then, we can prove:
\bt
Let $M(t)=\graph u(t)$ be the leaves of the inverse curvature flow, where $F$ and the initial hypersurface are smooth, then the estimate 
\begin{equation}
\norm{D^mA}\le c_m e^{-\frac tn}\qq\A\, m\ge 1
\end{equation}
is valid and the function
\begin{equation}\lae{1.12}
(u-2)e^{\frac tn}
\end{equation}
converges in $C^\un(\Ss[n])$ to a strictly negative function.
\et

\br
After publishing a first version of the paper in the arXiv we learnt that Qi Ding in \cite{ding:hyperbolic} published a similar result for the inverse mean curvature flow in $\Hh$ even claiming that the rescaled flow hypersurfaces would converge to a sphere. However, he used a somewhat crude rescaling, namely,
\begin{equation}\lae{1.13}
\frac{\breve u}t
\end{equation}
and not the finer
\begin{equation}\lae{1.14} 
\breve u-\frac tn 
\end{equation}
which we consider.

The fact that the functions in \re{1.13} converge to $\frac1n$ follows immediately from the estimate \fre{3.19}. 
\er

\section{Definitions and Conventions}

The main objective of this section is to state the equations of Gau{\ss}, Codazzi,
and Weingarten for hypersurfaces. For greater generality we shall formulate the governing equations of a
hypersurface $M$ in a semi-riemannian \di{(n+1)} manifold $N$, which is either
Riemannian or Lorentzian. Geometric quantities in $N$ will be denoted by
$(\bar g_{\alpha\beta}),(\riema \alpha\beta\gamma\delta)$, etc., and those in $M$ by $(g_{ij}), (\riem
ijkl)$, etc. Greek indices range from $0$ to $n$ and Latin from $1$ to $n$; the
summation convention is always used. Generic coordinate systems in $N$ resp.
$M$ will be denoted by $(x^\alpha)$ resp. $(\x^i)$. Covariant differentiation will
simply be indicated by indices, only in case of possible ambiguity they will be
preceded by a semicolon, i.e., for a function $u$ in $N$, $(u_\alpha)$ will be the
gradient and
$(u_{\alpha\beta})$ the Hessian, but e.g., the covariant derivative of the curvature
tensor will be abbreviated by $\riema \alpha\beta\gamma{\delta;\e}$. We also point out that
\begin{equation}
\riema \alpha\beta\gamma{\delta;i}=\riema \alpha\beta\gamma{\delta;\e}x_i^\e
\end{equation}
with obvious generalizations to other quantities.

Let $M$ be a \tit{spacelike} hypersurface, i.e., the induced metric is Riemannian,
with a differentiable normal $\n$. We define the signature of $\n$, $\sigma=\s(\n)$, by
\begin{equation}
\sigma=\bar g_{\alpha\beta}\n^\alpha\n^\beta=\spd \n\n.
\end{equation}
In case $N$ is Lorentzian, $\sigma=-1$, and $\n$ is time-like.

In local coordinates, $(x^\alpha)$ and $(\x^i)$, the geometric quantities of the
spacelike hypersurface $M$ are connected through the following equations
\begin{equation}\lae{2.3}
x_{ij}^\alpha=-\sigma h_{ij}\n^\alpha
\end{equation}
the so-called \tit{Gau{\ss} formula}. Here, and also in the sequel, a covariant
derivative is always a \tit{full} tensor, i.e.,

\begin{equation}
x_{ij}^\alpha=x_{,ij}^\alpha-\ch ijk x_k^\alpha+\cha \beta\gamma\alpha x_i^\beta x_j^\gamma.
\end{equation}
The comma indicates ordinary partial derivatives.

In this implicit definition the \tit{second fundamental form} $(h_{ij})$ is taken
with respect to $-\sigma\n$.

The second equation is the \tit{Weingarten equation}
\begin{equation}
\n_i^\alpha=h_i^k x_k^\alpha,
\end{equation}
where we remember that $\n_i^\alpha$ is a full tensor.

Finally, we have the \tit{Codazzi equation}
\begin{equation}
h_{ij;k}-h_{ik;j}=\riema\alpha\beta\gamma\delta\n^\alpha x_i^\beta x_j^\gamma x_k^\delta
\end{equation}
and the \tit{Gau{\ss} equation}
\begin{equation}
\riem ijkl=\sigma \{h_{ik}h_{jl}-h_{il}h_{jk}\} + \riema \alpha\beta\gamma\delta x_i^\alpha x_j^\beta x_k^\gamma
x_l^\delta.
\end{equation}
Here, the signature of $\n$ comes into play.

Now, let us assume that
$N$ is a topological product $\R[]\times \mc S_0$, where $\mc S_0$ is a
compact Riemannian manifold, and that there exists a Gaussian coordinate system
$(x^\alpha)$, such that the metric in $N$ has the form 
\begin{equation}\lae{2.8b}
d\bar s_N^2=e^{2\psi}\{\s{dx^0}^2+\sigma_{ij}(x^0,x)dx^idx^j\},
\end{equation}
where $\sigma_{ij}$ is a Riemannian metric, $\psi$ a function on $N$, and $x$ an
abbreviation for the spacelike components $(x^i)$,

We also assume that
the coordinate system is \tit{future oriented}, i.e., the time coordinate $x^0$
increases on future directed curves. Hence, the \tit{contravariant} time-like
vector $(\x^\alpha)=(1,0,\dotsc,0)$ is future directed as is its \tit{covariant}
version
$(\x_\alpha)=e^{2\psi}(\s,0,\dotsc,0)$.

Let $M=\graph \fv u\so$ be a spacelike hypersurface
\begin{equation}
M=\set{(x^0,x)}{x^0=u(x),\,x\in\mc S_0},
\end{equation}
then the induced metric has the form
\begin{equation}
g_{ij}=e^{2\psi}\{\s u_iu_j+\sigma_{ij}\}
\end{equation}
where $\sigma_{ij}$ is evaluated at $(u,x)$, and its inverse $(g^{ij})=(g_{ij})^{-1}$ can
be expressed as
\begin{equation}\lae{2.10}
g^{ij}=e^{-2\psi}\{\sigma^{ij}-\s \frac{u^i}{v}\frac{u^j}{v}\},
\end{equation}
where $(\sigma^{ij})=(\sigma_{ij})^{-1}$ and
\begin{equation}\lae{2.11}
\begin{aligned}
u^i&=\sigma^{ij}u_j\\
v^2&=1+\s \sigma^{ij}u_iu_j\equiv 1+\s \abs{Du}^2.
\end{aligned}
\end{equation}

The covariant form of a normal vector of a graph looks like
\begin{equation}
(\n_\alpha)=\pm v^{-1}e^{\psi}(1, -u_i).
\end{equation}
and the contravariant version is
\begin{equation}
(\n^\alpha)=\pm v^{-1}e^{-\psi}(\s, -u^i).
\end{equation}

In the Gau{\ss} formula \re{2.3} we are free to choose any of two normals, but we stipulate that in general we use \begin{equation}
(\n^\alpha)= v^{-1}e^{-\psi}(\s, -u^i).
\end{equation}
as normal vector. 

Look at the component $\alpha=0$ in \re{2.3}, then we obtain 
\begin{equation}\lae{2.16}
e^{-\psi}v^{-1}h_{ij}=-u_{ij}-\cha 000\mspace{1mu}u_iu_j-\cha 0i0
\mspace{1mu}u_j-\cha 0j0\mspace{1mu}u_i-\cha ij0.
\end{equation}
Here, the covariant derivatives a taken with respect to the induced metric of
$M$, and
\begin{equation}\lae{2.18}
-\cha ij0=e^{-\psi}\bar h_{ij},
\end{equation}
where $(\bar h_{ij})$ is the second fundamental form of the hypersurfaces
$\{x^0=\const\}$.

\section{First Estimates}

Let $F\in C^{m,\al}(\C)$, $m\ge 4$, be a monotone and concave curvature function, homogeneous of degree $1$, and normalized such that
\begin{equation}\lae{3.1}
F(1,\ldots,1)=n.
\end{equation}

We first look at the flow of a geodesic sphere $S_{r_0}$. Fix a point $p_0\in \Hh$ and consider geodesic polar coordinates centered at $p_0$. Then the hyperbolic metric can be expressed as
\begin{equation}\lae{3.2}
d\bar s^2=dr^2+\sinh^2r\s_{ij}\,dx^idx^j,
\end{equation}
where $\s_{ij}$ is the canonical metric of $\Ss[n]$.

Geodesic spheres $S_r$ with center in $p_0$ are umbilic and their second fundamental form is given by
\begin{equation}
\bar h_{ij}=\coth r\bar g_{ij},
\end{equation}
where
\begin{equation}
\bar g_{ij}=\sinh^2r\s_{ij}.
\end{equation}
Hence, if we consider an inverse curvature flow (ICF) with initial hypersurface $S_{r_0}$, then the flow hypersurfaces $M(t)$ will be spheres with radii $r(t)$ satisfying the scalar curvature flow equation
\begin{equation}
\dot r=\frac1F=\frac1{n\coth r},
\end{equation}
and we deduce further, from
\begin{equation}\lae{3.6}
\coth r\, dr=\frac1n dt,
\end{equation}
\begin{equation}
\log\sinh r-\log\sinh r_0=\frac tn,
\end{equation} 
or equivalently,
\begin{equation}\lae{3.8}
\sinh r=\sinh r_0 e^\frac{t}n.
\end{equation}
  
Let us now consider the inverse curvature flow of a star-shaped hypersurface $M_0$ which is given as a graph over $\Ss[n]$
\begin{equation}
M_0=\fv{\graph u_0}{\Ss[n]}.
\end{equation}
The flow exists on a maximal time interval $[0,T^*)$, $0<T^*\le\un$, and its leaves are also graphs
\begin{equation}
M(t)=\fv{\graph u(t)}{\Ss[n]},
\end{equation}
which satisfy, besides the original flow equation,
\begin{equation}
\dot x=-\F\nu=\frac1F\nu,
\end{equation}
the scalar flow equation
\begin{equation}
\dot u=\frac{\tilde v}F,
\end{equation}
where
\begin{equation}
\tilde v=v^{-1},\q v^2=1+\abs{Du}^2=1+\frac1{\sinh^2u}\s^{ij}u_iu_j,
\end{equation}
where the dot indicates a total time derivative. If we instead consider a partial time derivate, then we get
\begin{equation}\lae{3.14} 
\dot u\equiv\pde ut=\frac vF, 
\end{equation}
\cf \cite[p. 98]{cg:cp}.

Let $S_{r_i}$, $i=1,2$, be geodesic spheres satisfying
\begin{equation}\lae{3.15}
r_1<u_0<r_2,
\end{equation}
and let $u_i$, $i=1,2$, be the solutions to the corresponding inverse curvature flows, then this inequality will also be valid for $t>0$, i.e.,
\begin{equation}
u_1(t)<u(t)<u_2(t)\qq\A\,t\in [0,T^*),
\end{equation}
in view of the maximum principle, and we conclude:
\bl\lal{3.1}
The solutions $M(t)=\graph u(t)$ of the ICF satisfy the estimates
\begin{equation}\lae{3.17}
\sinh r_1<\sinh u(t) e^{-\frac tn}<\sinh r_2\qq\A\, t\in [0,T^*),
\end{equation}
and there exist constants $c_i$, $i=1,2$, such that the function
\begin{equation}\lae{3.18}
\tilde u=u-\frac tn
\end{equation}
is uniformly bounded by
\begin{equation}\lae{3.19}
c_1<\tilde u(t)<c_2\qq\A\,t\in [0,T^*).
\end{equation}
\el

\bp
The inequality \re{3.17} follows from \re{3.8} and the parabolic maximum principle, while \re{3.19} is due to the trivial estimate
\begin{equation}
0<\tilde c_1\le \sinh r \,e^{-r}\le \tilde c_2\qq\A\, 0<r_0\le r
\end{equation}
with appropriate constants $\tilde c_i$.
\ep

Next, we want to derive an a priori estimate for $v$, or equivalently, for
\begin{equation}\lae{3.21}
\abs{Du}^2=\frac1{\sinh^2u}\s^{ij}u_iu_j.
\end{equation}

Let us write the metric \re{3.2} in a more general form
\begin{equation}\lae{3.22}
d\bar s^2=dr^2+\vt^2(r)\s_{ij} dx^idx^j.
\end{equation}

The second fundamental form of $\graph u$ can then be expressed as
\begin{equation}
\begin{aligned}
h_{ij}v^{-1}&=-u_{ij}+\bar h_{ij}\\
&=-u_{ij}+\dot\vt\vt\s_{ij}.
\end{aligned}
\end{equation}

Define the metric
\begin{equation}
\tilde \s_{ij}=\vt^2(u)\s_{ij},
\end{equation}
and denote covariant differentiation with respect to this metric by a semi-colon, then
\begin{equation}
h_{ij}v^{-1}=-v^{-2}u_{;ij}+\dot\vt\vt \s_{ij},
\end{equation}
\cf \cite[Lemma 2.7.6]{cg:cp}, and we conclude further
\begin{equation}\lae{3.26}
\begin{aligned}
h^i_j&=g^{ik}h_{kj}\\
&=v^{-1}\vt^{-1}\{-(\s^{ik}-v^{-2}\f^i\f^k)\f_{jk}+\dot\vt \de^i_j\},
\end{aligned}
\end{equation}
where $\s^{ij}$ is the inverse of $\s_{ij}$,
\begin{equation}
\f=\int_{r_0}^u\vt^{-1},
\end{equation}
\begin{equation}
\f^i=\s^{ik}\f_k,
\end{equation}
and $\f_{jk}$ are the second covariant derivatives of $\f$ with respect to the metric $\s_{ij}$.

Thus, the scalar curvature equation \re{3.14} can now be expressed as
\begin{equation}
\dot u=\frac v{F(h^i_j)},
\end{equation}
or equivalently,
\begin{equation}\lae{3.30} 
\dot\f=\vt^{-1}\dot u =\frac1{F(\vt v^{-1}h^i_j)}\equiv \frac1{F(\tilde h^i_j)},
\end{equation}
where
\begin{equation}
\tilde h^i_j=v^{-2}\{-(\s^{ik}-v^{-2}\f^i\f^k)\f_{jk}+\dot\vt \de^i_j\}.
\end{equation}
Let
\begin{equation}
\tilde g_{ij}=\f_i\f_j+\s_{ij},
\end{equation}
then we consider the eigenvalues of
\begin{equation}
\tilde h_{ij}=\tilde g_{ik}\tilde h^k_j
\end{equation}
with respect to this metric and we define $F^{ij}$ \resp $F^i_j$ accordingly
\begin{equation}
F^{ij}=\pdc F{\tilde h}{{ij}}
\end{equation}
and
\begin{equation}
F^i_j=\pdm F{\tilde h}ij=\tilde g_{jk}F^{ik}.
\end{equation}
Note that $\tilde h_{ij}$ is symmetric, since $h_{ij}$ and $\tilde g_{ij}$ can be diagonalized simultaneously. We also emphasize that
\begin{equation}\lae{3.36}
\abs{Du}^2=\s^{ij}\f_i\f_j\equiv\abs{D\f}^2.
\end{equation}

\bl\lal{3.2}
Let $u$ be solution of the scalar curvature equation
\begin{equation}
\dot u=\frac v{F(h_{ij})},
\end{equation}
then
\begin{equation}\lae{3.38}
\abs{Du}^2\le c.
\end{equation}
Moreover, if $F$ is bounded from above
\begin{equation}
F\le c_0,
\end{equation}
then there exists $0<\lam=\lam(c_0)$ such that
\begin{equation}\lae{3.40}
\abs{Du}^2\le c e^{-\lam t}\qq\A\, t\in [0,T^*).
\end{equation}
\el
\bp
\cq{\re{3.38}}\q In view of  \re{3.36}, we may estimate 
\begin{equation}
w=\tfrac12\abs{D\f}^2.
\end{equation}
Differentiating equation \re{3.30} covariantly with respect to
\begin{equation}
\f^kD_k
\end{equation}
we deduce 
\begin{equation}\lae{3.43}
\begin{aligned}
\dot w&=F^{-2}\{2v^{-2}Fw_i\f^i+v^{-2}F^k_l\tilde g^{lr}w_{kr} -v^{-2}F^k_l\tilde g^{lr}\f_{ik}\f^i_r\\
&\q+v^{-2}F^k_l\tilde g^{lr}_{\hp{lr};i}\f^i\f_{kr}+v^{-2}F^k_l\tilde g^{lr}\f_r\f_k-v^{-2}F^k_l\tilde g^{lr}\s_{kr}\abs{D\f}^2\\
&\q -2v^{-2}F^k_k\Ddot\vt\vt w\},
\end{aligned}
\end{equation}
where covariant derivatives with respect to the metric $\s_{ij}$ are simply denoted by indices, if no ambiguities are possible, and by a semi-colon otherwise. In deriving the previous equation we also used the Ricci identities and the properties of the Riemann curvature tensor of $\Ss[n]$.

Now, let $0<T<T^*$ and suppose that
\begin{equation}
\sup_{Q_T}w,\qq Q_T=[0,T]\times \Ss[n],
\end{equation}
is attained at $(t_0,x_0)$ with $t_0>0$. Then the maximum principle implies
\begin{equation}\lae{3.45}
\begin{aligned}
0&\le v^{-2}\{-F^k_l\tilde g^{lr}\f_{ik}\f^i_r+(F^k_l\tilde g^{lr}\f_r\f_k-F^k_l\tilde g^{lr}\s_{kr}\abs{D\f}^2)\\
&\q -2F^k_k\sinh^2u\, w\}.
\end{aligned}
\end{equation}
The right-hand side, however, is strictly negative, if $w>0$, hence $t_0>0$ is not possible, since we didn't assume $M_0$ to be a sphere, and we conclude
\begin{equation}
w\le \sup_{\Ss[n]}w(0).
\end{equation}

\cq{\re{3.40}}\q Now, assume that the original curvature function is uniformly bounded
\begin{equation}
F(h_{ij})\le c_0,
\end{equation}
and let $0<\lam$ be a constant, then
\begin{equation}
\tilde w=we^{\lam t}
\end{equation}
satisfies the same equation as $w$ with an additional term
\begin{equation}
\lam \tilde w
\end{equation}
at the right-hand side.

Applying the maximum principle as before, we deduce, that at a point $(t_0,x_0)$, $t_0>0$, where $\tilde w$ attains a positive maximum, there holds instead of \re{3.45}
\begin{equation}
0< -2F^k_k\sinh^2u\,\tilde w+\lam v^2F^2(\tilde h^i_j)\tilde w,
\end{equation}
but
\begin{equation}
vF(\tilde h^i_j)=\sinh uF(h_{ij})\le \sinh u\,c_0,
\end{equation}
and hence
\begin{equation}
we^{\lam t}\le \sup_{\Ss[n]}w(0)
\end{equation}
for all 
\begin{equation}
0<\lam\le 2nc^{-2}_0,
\end{equation}
since
\begin{equation}
F^k_k\ge n.
\end{equation}
\ep

\section{$C^2$-estimates and existence for all time}

To prove estimates for $h_{ij}$, we first need an a priori bound for $F$.
\bl
Let $M(t)$ be the leaves of the ICF
\begin{equation}
\dot x=-\F\nu,
\end{equation}
then there exists a positive constant $c_1$ such that
\begin{equation}\lae{4.2}
0<c_1\le F\qq\A\, t\in [0,T^*).
\end{equation}
\el
\bp
The function $\F$, or equivalently $-\F$, satisfies the linear parabolic equation
\begin{equation}\lae{4.3}
\begin{aligned}
\F'-\dot\F F^{ij}\F_{ij}=\dot\F F^{ij}h_{ik}h^k_j\F +K_N\dot\F F^{ij}g_{ij}\F,
\end{aligned}
\end{equation}
when the ambient Riemannian space $N$ is a space of constant curvature $K_N$, \cf \cite[Corollary 3.5]{cg:spaceform}.

Another very useful equation is satisfied by a quantity $\chi$ which is defined by
\begin{equation}
\chi=v\h,
\end{equation}
where $0<\h=\h(r)$ is a solution of
\begin{equation}
\dot\h=-\frac{\bar H}n\h;
\end{equation}
here $r$ is the radial distance to  the center of geodesic polar coordinates in a spaceform $N$, and $\bar H=\bar H(r)$ is the mean curvature of $S_r$.

When $N=\Hh$, $\h$ is given by
\begin{equation}
\h=\frac1{\sinh r},
\end{equation}
and 
\begin{equation}
\chi = v\h(u)
\end{equation}
then satisfies
\begin{equation}
\begin{aligned}
\dot\chi -\dot\F F^{ij}\cchi{ij}=-\dot\F F^{ij}h_{ik}h^k_j-2\chi^{-1}\dot\F F^{ij}\cchi i\cchi j+\{\dot\F F+\F\}\frac{\bar H}nv\chi
\end{aligned}
\end{equation}
for a general function $\F$. In case of the inverse curvature, the term in the braces on the right-hand side vanishes.

In view of \rl{3.1}, the function
\begin{equation}\lae{4.9}
\tilde\chi=\chi e^\frac tn
\end{equation}
is uniformly bounded,
\begin{equation}
0<c_1\le \tilde\chi \le c_2\qq\A\, t\in[0,T^*),
\end{equation}
and
\begin{equation}\lae{4.11}
\dot{\tilde\chi}-\dot\F F^{ij}\tchi{ij}=- \dot\F F^{ij} h_{ik}h^k_j\tilde \chi-2\tilde\chi^{-1} \dot\F F^{ij}\tchi i\tchi j+\tfrac1n \tilde\chi,
\end{equation}
when we consider an ICF.

We claim that
\begin{equation}\lae{4.12}
w=\log(-\F)+\log\tilde\chi\le\const
\end{equation}
during the evolution, which in turn would prove \re{4.2}.

To derive \re{4.12} we first fix $0<T<T^*$ and let 
\begin{equation}
(t_0,\xi_0)\in Q_T=[0,T]\times \Ss[n],\q t_0>0,
\end{equation}
be such that
\begin{equation}
w(t_0,\xi_0)=\sup_{Q_T}w.
\end{equation}

The equations \re{4.3}, \re{4.11} and the maximum principle then yield in $(t_0,\xi_0)$
\begin{equation}
0\le -\dot\F F^{ij}g_{ij}+\tfrac1n,
\end{equation}
which can only hold, if
\begin{equation}
n\le F(t_0,\xi_0);
\end{equation}
hence $w$ is uniformly bounded from above.
\ep

\bl
During the evolution $F$ is uniformly bounded from above.
\el
\bp
The function $u$ satisfies the parabolic equation
\begin{equation}
\begin{aligned}
\dot u-\dot\F F^{ij}u_{ij}&=-\F v^{-1}+\dot\F F v^{-1}-\dot\F F^{ij}\bar h_{ij}\\
&= 2\dot\F F v^{-1} -\dot\F F^{ij}\bar h_{ij},
\end{aligned}
\end{equation}
\cf \cite[Lemma 3.3.2]{cg:cp}, and the rescaled function
\begin{equation}
\tilde u=u-\frac tn
\end{equation}
is uniformly bounded, \cf \rl{3.1}, and there holds
\begin{equation}\lae{4.19}
\dot{\tilde u}-\dot\F F^{ij}\tilde u_{ij}=2\dot\F F v^{-1}-\dot\F F^{ij}\bar h_{ij}-\tfrac1n.
\end{equation}

The lemma will be proved, if we can show
\begin{equation}
w=-\log(-\F)+\tilde u=\log F+\tilde u\le\const
\end{equation}
during the evolution.

Applying the maximum principle as before, we conclude
\begin{equation}
0\le 2F^{-1}v^{-1}-\tfrac1n,
\end{equation}
hence, $F$ has to be bounded, proving the claim.
\ep

As an immediate corollary we deduce, in view of \rl{3.2}: 
\br
$\abs{Du}^2$ satisfies the estimate \re{3.40}, i.e., it decays exponentially, if $T^*=\un$. 
\er

We are now ready to prove a priori estimates for the principal curvatures $\ka_i$. The proof will be similar to a corresponding proof in \cite[Theorem 1.4]{cg:weingarten} valid in arbitrary Riemannian spaces. Our former result cannot be applied directly, since we assumed that the flow stays in a compact subset and also considered a contracting flow not an expanding one as we do now.

\bl\lal{4.4}
The principal curvatures of the flow hypersurfaces are uniformly bounded from above
\begin{equation}
\ka_i\le \const\qq\A\, 1\le i\le n,
\end{equation}
and hence, are compactly contained in $\C,$ in view of the estimate \re{4.2}.
\el
\bp
In a Riemannian space of constant curvature the second fundamental forms $h^i_j$ of the flow hypersurfaces $M(t)$ satisfy the evolution equation
\begin{equation}\lae{4.23}
\begin{aligned}
\dot h^i_j-\dot\F F^{kl}h^i_{j;kl}&=\dot\F F^{kl}h_{kr}h^r_kh^i_j+(\F-\dot\F F)h^k_ih_{kj}+\Ddot\F F_jF^i\\
&\q \dot\F F^{kl,rs}h_{kl;i}h_{rs;}^{\hp{rs;}i}+K_N\{(\F+\dot\F F)\de^i_j-\dot\F F^{kl}g_{kl}h^i_j\}
\end{aligned}
\end{equation}
\cf \cite[Lemma 2.4.3]{cg:cp}.

Here, the flow is given as an embedding
\begin{equation}
x=x(t,\xi),\qq (t,\xi)\in [0,T^*)\ti \Ss[n],
\end{equation}
and
\begin{equation}
F_i=\pd F{\xi}i=F^{kl}h_{kl;i}.
\end{equation}

By assumption, $F$ is monotone and concave. Thus, choosing, in a given point, coordinates $(\xi^i)$ such that
\begin{equation}
g_{ij}=\de_{ij}\q\wed\q h_{ij}=\ka_i\de_{ij},
\end{equation}
and labelling the $\ka_i$ such that
\begin{equation}\lae{4.27}
\ka_1\le\cdots\le \ka_n,
\end{equation}
then
\begin{equation}\lae{4.28}
\begin{aligned}
F^{kl,rs}\h_{kl}\h_{rs}&\le \sum_{k\not=l}\frac{F^{kk}-F^{ll}}{\ka_k-\ka_l}(\h_{kl})^2\\
&\le \frac2{\ka_n-\ka_1}\sum_{k=1}^n(F^{nn}-F^{kk})(\h_{nk})^2,
\end{aligned}
\end{equation}
and
\begin{equation}\lae{4.29}
F^{nn}\le\cdots\le F^{11}.
\end{equation}
For a proof of \re{4.28} see \cite[Lemma 1.1]{cg:weingarten} and of \re{4.29} \cite[Lemma 2]{eh2}. 

Let $\tilde\chi$ be the rescaled function in \re{4.9} and define 
\begin{equation}
\hat\chi =\tilde\chi^{-1},
\end{equation}
then there exists a constant $\theta>0$ such that
\begin{equation}
2\theta \le \hat\chi.
\end{equation}

Next, let $\zeta,\f$ and $w$ be defined by
\begin{equation}\lae{4.32}
\zeta=\sup\set{h_{ij}\h^i\h^j}{\norm\h=1},
\end{equation}
\begin{equation}
\f=-\log(\hat\chi-\theta)
\end{equation}
and 
\begin{equation}
w=\log\zeta+\f+\lam\tilde u,
\end{equation}
where $\tilde u$ is the function in \fre{3.18} and $\lam>0$ is supposed to be large.  We claim that $w$ is bounded, if $\lam$ is chosen sufficiently large.

Let $0<T<T^*$, and $x_0=x_0(t_0,\xi_0)$, with $ 0<t_0\le T$, be a point in $M(t_0)$ such
that
\begin{equation}
\sup_{M_0}w<\sup\set {\sup_{M(t)} w}{0<t\le T}=w(x_0).
\end{equation}

We then introduce a Riemannian normal coordinate system $(\x^i)$ at $x_0\in
M(t_0)$ such that at $x_0=x(t_0,\x_0)$ we have
\begin{equation}
g_{ij}=\delta_{ij}\q \tup{and}\q \zeta=h_n^n.
\end{equation}

Let $\tilde \h=(\tilde \h^i)$ be the contravariant vector field defined by 
\begin{equation}
\tilde \h=(0,\dotsc,0,1),
\end{equation}
and set
\begin{equation}
\tilde \zeta=\frac{h_{ij}\tilde \h^i\tilde \h^j}{g_{ij}\tilde \h^i\tilde \h^j}\raise 2pt
\hbox{.}
\end{equation}

$\tilde\zeta$ is well defined in neighbourhood of $(t_0,\x_0)$.

Now, define $\tilde w$ by replacing $\zeta$ by $\tilde \zeta$ in \re{4.32}; then, $\tilde w$
assumes its maximum at $(t_0,\x_0)$. Moreover, at $(t_0,\x_0)$ we have 
\begin{equation}
\dot{\tilde \zeta}=\dot h_n^n,
\end{equation}
and the spatial derivatives do also coincide; in short, at $(t_0,\x_0)$ $\tilde \zeta$
satisfies the same differential equation \re{4.23} as $h_n^n$. For the sake of
greater clarity, let us therefore treat $h_n^n$ like a scalar and pretend that $w$
is defined by 
\begin{equation}
w=\log h_n^n+ \f+\lam\tilde u.
\end{equation} 

From equations \re{4.23}, \re{4.28}, \re{4.11} and \re{4.19} we infer that in $(t_0,\xi_0)$ 
\begin{equation}\lae{4.41}
\begin{aligned}
0&\le -\dot\F F^{ij}h_{ik}h^k_j\frac\theta{\hat\chi-\theta}+(\F-\dot\F F)h^n_n-(\F+\dot\F F)(h^n_n)^{-1}\\
&\q+\dot\F F F^{kl}g_{kl}   +(\F+\dot\F F)\frac{\bar H}nv\frac{\hat\chi}{\hat\chi-\theta}+\frac1n\frac{\hat\chi}{\hat\chi-\theta}\\
&\q +\lam(-\F+\dot\F F)v^{-1}-\lam \dot\F F^{ij}\bar h_{ij}-\frac\lam{n} \\
&\q +\dot\F F^{ij}(\log h^n_n)_i(\log h^n_n)_j-\dot\F F^{ij}\f_i\f_j\\
&\q +\frac2{\ka_n-\ka_1}\dot\F\sum_{i=1}^n(F^{nn}-F^{ii})(h^{\hp{ni;}n}_{ni;})^2(h^n_n)^{-1}.
\end{aligned}
\end{equation}

There holds
\begin{equation}
F^{ij}\bar h_{ij}\ge c_0F^{ij}g_{ij},\q c_0>0;
\end{equation}
moreover,
\begin{equation}
h_{ni;n}=h_{nn;i},
\end{equation}
and
\begin{equation}
\F+\dot\F F=0,
\end{equation}
though
\begin{equation}
\F\le 0\q\wed\q \abs{\F}\le c\dot\F F
\end{equation}
would suffice.

We then distinguish two cases.

\cvm
\tit{Case} $1$.\q Suppose that
\begin{equation}
\ka_1< -\e_1 \ka_n,
\end{equation}
where $\e_1>0$ is small, note that the principal curvatures are labelled according to \re{4.27}. Then, we infer from \cite[Lemma 8.3]{cg:scalar}
\begin{equation}
F^{ij}h_{ki}h^k_j\ge \tfrac1n F^{ij}g_{ij}\e_1^2\ka_n^2,
\end{equation}
and 
\begin{equation}
F^{ij}g_{ij}\ge F(1,\ldots,1),
\end{equation}
for a proof see e.g., \cite[Lemma 2.2.19]{cg:cp}.

Since $Dw=0$,
\begin{equation}
D\log h^n_n=-D\f-\lam D\tilde u,
\end{equation}
we obtain
\begin{equation}
\dot\F F^{ij}(\log h^n_n)_i(\log h^n_n)_j=\dot \F F^{ij}\f_i\f_j+2\lam \dot\F F^{ij}\f_i\tilde u_j+\lam^2\dot\F F^{ij}\tilde u_i\tilde u_j,
\end{equation}
where
\begin{equation}
\abs{\f_i}\le c\abs{\ka_i}\norm{Du}+c\norm{Du},
\end{equation}
as one easily checks. 

Hence, we conclude that $\ka_n$ is a priori bounded in this case for any choice of $\lam >0$, if we use
\begin{equation}
F\le \const,
\end{equation}
or for $\lam >2$ otherwise.

Let us remark that
\begin{equation}
\frac{\hat\chi}{\hat\chi-\theta}\le 2
\end{equation}
and
\begin{equation}
F\le F(1,\dots,1)\ka_n=n\ka_n.
\end{equation}

\cvm
\tit{Case} 2. Suppose that
\begin{equation}
\ka_1\ge -\e_1\ka_n,
\end{equation}
then the last term in inequality \re{4.41} can be estimated from above by
\begin{equation}
\frac2{1+\e_1}\dot\F \sum_{i=1}^n(F^{nn}-F^{ii})(\log h^n_{n;i})^2.
\end{equation}

The terms in \re{4.41} containing derivatives of $h^n_n$   can therefore be estimated from above by
\begin{equation}
\begin{aligned}
&-\frac{1-\e_1}{1+\e_1}\dot\F F^{ij}(\log h^n_n)_i(\log h^n_n)_j+\frac1{1+\e_1}\dot\F F^{nn}\sum_{i=1}^n(\log h^n_{n;i})^2\\
&\qq\qq\le \dot\F F^{nn}\sum_{i=1}^n(\log h^n_{n;i})^2\\
&\qq\qq=\dot\F F^{nn}\norm{D\f+\lam D\tilde u}^2\\
&\qq\qq=\dot\F F^{nn}\{\norm{D\f}^2+\lam^2\norm{D\tilde u}^2+2\lam \spd{D\f}{D\tilde u}\}.
\end{aligned}
\end{equation}

Hence, we finally deduce
\begin{equation}
\begin{aligned}
0&\le -\dot\F F^{nn}\ka_n^2\frac\theta{\hat\chi-\theta}-\dot\F F\ka_n+\dot\F F^{kl}g_{kl}(1-\lam c_0)+c\\
&\q+\lam c\dot\F F -\frac\lam{n}+\lam^2 c\dot\F F^{nn}(1+\ka_n).
\end{aligned}
\end{equation}
Thus, we obtain an a priori estimate
\begin{equation}
\ka_n\le\const,
\end{equation}
if $\lam$ is chosen large enough. Note that $\e_1$ is only subject to the requirement
\begin{equation}
0<\e_1<1.
\end{equation} 
\ep

As a corollary we can state:
\bc
Let the initial hypersurface $M_0\in C^{m+2,\al}$, $4\le m\le\un$, $0<\al<1$, then the solution of the curvature flow
\begin{equation}
\dot x=-\F\nu
\end{equation}
exists for all time and belongs to the parabolic H\"older space 
\begin{equation}
H^{m+\al,\frac{m+\al}2}(Q),
\end{equation}
while the solution $u$ of the scalar flow belongs to
\begin{equation}
H^{m+2+\al,\frac{m+2+\al}2}(Q),
\end{equation} 
where
\begin{equation}
Q=[0,\un)\times \Ss[n].
\end{equation}
The norm will still depend on $t$ however due to the present coordinate system.
\ec

\bp
Let us look at the scalar flow equation \fre{3.14}. In view of the previous estimates the nonlinear operator is uniformly elliptic and by assumption also concave, hence we may apply the Krylov-Safonov estimates yielding uniform H\"older estimates for $\dot u$ and $D^2u$ estimates. Now, the linear theory and the parabolic Schauder estimates can be applied; for details see  e.g.\ \cite[Chapter 2.6]{cg:cp} and \cite[Section 6]{cg:survey}.
\ep
\section{The conformally flat parametrization}
Hyperbolic space is conformally flat such that
\begin{equation}
\begin{aligned}
d\bar s^2&=\frac1{(1-\frac14\abs x^2)^2}dx^2\\
&=\frac1{(1-\frac14 r^2)^2}\{dr^2+r^2\s_{ij}dx^idx^j\}\\
&\equiv e^{2\psi}\{dr^2+r^2\s_{ij}dx^idx^j\}
\end{aligned}
\end{equation}
after introducing polar coordinates.

Define the variable $\tau$ by
\begin{equation}
d\tau=\frac1{1-\frac14 r^2}dr
\end{equation}
such that
\begin{equation}
\tau=\log(2+r)-\log(2-r),
\end{equation}
then 
\begin{equation}\lae{5.4} 
\sinh^2\tau=\frac{r^2}{(1-\frac14 r^2)^2},
\end{equation}
and we see that $\tau$ is the radial distance in hyperbolic space from the origin of the euclidean ball $B_2(0)$.

A star-shaped hypersurface $M\su\Hh$ is also star-shaped in $\R$ under this correspondence.

Let us distinguish geometric quantities in $\Hh$ by an additional breve from the corresponding quantities in $\R$, e.g., $\breve g_{\al\bet}$, $\breve g_{ij}$, $M=\graph \breve u$, $\breve h_{ij}$, $\breve\nu$, etc.

Consider a hypersurface
\begin{equation}
M=\graph\breve u=\graph u,
\end{equation}
then
\begin{equation}
\breve u=\log(2+u)-\log(2-u)
\end{equation}
and
\begin{equation}
\breve u_i=\frac1{1-\frac14 u^2} u_i
\end{equation}
and $\abs{D\breve u}^2$ as defined in \fre{3.21} can be expressed as
\begin{equation}
\abs{D\breve u}^2=u^{-2}\s^{ij}u_iu_j\equiv\abs{Du}^2,
\end{equation}
hence the term $v$ is identical in both coordinate systems which is also evident from the invariant definition of $v$ by
\begin{equation}
v^{-1}=\spd\h{\nu},
\end{equation}
where
\begin{equation}
\h=Dd
\end{equation}
and $d$ is the distance function in hyperbolic space from the origin.

The second fundamental forms are connected through the relation
\begin{equation}\lae{5.11}
\begin{aligned}
e^\psi\breve h^i_j&=h^i_j+\psi_\al\nu^\al\de^i_j\\
&\equiv h^i_j+v^{-1}\tilde\vt\de^i_j,
\end{aligned}
\end{equation}
where
\begin{equation}
\tilde\vt=\tfrac12\fr{r}r.
\end{equation}

Let
\begin{equation}\lae{5.13}
\check h_{ij}=h_{ij}+v^{-1}\tilde\vt g_{ij},
\end{equation}
\begin{equation}
g_{ij}=u_iu_j+u^2\s_{ij},
\end{equation}
then the curvature flow in $\Hh$
\begin{equation}
\dot x=F^{-1}\breve\nu
\end{equation}
can also be viewed as a curvature flow in $\R$
\begin{equation}\lae{5.16}
\dot x=F^{-1}\nu,
\end{equation}
where now $F$ depends on the eigenvalues of $\check h_{ij}$ with respect to the metric $g_{ij}$
\begin{equation}
F=F(\check h_{ij})=F(\check h^i_j).
\end{equation}

For the rest of this paper we shall mainly consider the curvature flow \re{5.16}.

Let us quickly summarize the most important flow equations.

Writing \re{5.16} slightly more general
\begin{equation}
\dot x=-\F(F)\nu\equiv-\F\nu
\end{equation}
there holds
\begin{equation}
\dot h^j_i=\F^j_i+\F h^k_i h^j_k,
\end{equation}
\cf \cite[Lemma 2.3.3]{cg:cp}, which will be the main ingredient to derive the subsequent modified flow equations:
\begin{equation}
\begin{aligned}
\F'-\dot\F F^{ij}\F_{ij}&= \dot\F F^{ij}h_{ki}h^k_j\F-\dot\F F^{ij}g_{ij}r_{\al\bet}\nu^\al\nu^\bet \tilde\vt\F \\
&\q -\dot\F F^{ij}g_{ij}\dot{\tilde\vt}v^{-2}\F+\dot\F F^{ij}g_{ij} \tilde\vt \F_ku^k,
\end{aligned}
\end{equation}
\begin{equation}\lae{5.21}
\begin{aligned}
\dot u-\dot\F F^{ij}u_{ij}=2v^{-1}F^{-1}-\dot\F F^{ij}g_{ij}\tilde\vt v^{-2}-\dot\F F^{ij}\bar h_{ij},
\end{aligned}
\end{equation}
where we used that 
\begin{equation}
\F(t)=-t^{-1},\qq t>0,
\end{equation}
here $t$ is just a symbol for a real variable, and where 
\begin{equation}
\bar h_{ij}=u^{-1}\bar g_{ij}=u\s_{ij}
\end{equation}
is the second fundamental form of the slices $\{x^0=u\}$, i.e., of spheres in $\R$ with center in the origin and radius $r=u$.

The evolution equation for the second fundamental form looks like
\begin{equation}\lae{5.24}
\begin{aligned}
\dot h^j_i-\dot\F F^{kl}h^j_{i;kl}&=\dot\F F^{kl}h_{kr}h^r_lh^j_i+(\F-\dot\F F)h_{ik}h^{kj}\\
&\q+\Ddot\F F_iF^j+\dot\F F^{kl,rs}\check h_{kl;i}\check h_{rs;}^{\hp{kl;}j}\\
&\q+ \dot\F F^{kl}g_{kl}\{-\tilde \vt v^{-1} h^r_ih_r^j-v^{-2}\dot{\tilde\vt}h^j_i+v^{-1}\dot{\tilde\vt}\bar h_{ik}g^{kj}\\
&\q-r_{\al\bet}\nu^\al\nu^\bet\tilde\vt h^j_i+\tilde\vt u_rh_{i;}^{r\hp{i;}j}+\dot{\tilde\vt}h^r_iu_ru^j+\dot{\tilde\vt}h^{rj}u_ru_i\\
&\q +r_{\al\bet}x^\al_kx^{\bet j}h^k_i\tilde\vt+v^{-1}\Ddot{\tilde\vt}u_iu^j+r_{\al\bet\ga} \nu^\al x^\bet_i x^\ga_j\tilde\vt\\
&\q +r_{\al\bet}\nu^\al x^\bet_i\dot{\tilde\vt}u^j+r_{\al\bet}\nu^\al x^{\bet j}\dot{\tilde\vt}u_i\}.
\end{aligned}
\end{equation}

The function $\chi$ is now defined by
\begin{equation}
\chi=v u^{-1}
\end{equation}
and there holds:
\bl
$\chi$ satisfies the evolution equation
\begin{equation}\lae{5.26}
\begin{aligned}
\dot\chi-\dot\F F^{ij}\chi_{ij}&=-\dot\F F^{ij}h^k_ih_{kj}\chi-2\chi^{-1}\dot\F F^{ij}\chi_i\chi_j+\{\dot\F F+\F\}\chi^2\\
&\q+\dot\F F^{ij}g_{ij}\{-\tilde\vt \chi^2+\chi_ku^ku\theta -\chi \dot\theta \norm{Du}^2u\},
\end{aligned}
\end{equation}
where
\begin{equation}
\theta(r)=r^{-1}\tilde\vt=\tfrac12\fr1r.
\end{equation}
\el
\bp
We consider a general $\F$ in the curvature flow in $\R$
\begin{equation}
\dot x=-\F\nu,
\end{equation}
where
\begin{equation}
F=F(\check h_{ij})
\end{equation}
and $\check h_{ij}$ is defined by \re{5.13}.

For the above flow the normal evolves according to
\begin{equation}
\dot\nu=\F^k x_k,
\end{equation}
\cf \cite[Lemma 2.3.2]{cg:cp}.

Using an euclidean coordinate system $(x^\al)$ in $\R$ it follows immediately that $\chi$ can be expressed as
\begin{equation}
\chi=\spd x\nu^{-1},
\end{equation}
and hence
\begin{equation}\lae{5.32}
\begin{aligned}
\dot\chi&=-\chi^2\spd{\dot x}\nu-\chi^2\spd x{\dot\nu}\\
&=\F\chi^2-\chi^2\F^ku_ku\\
&=\F\chi^2-\chi^2\dot\F F^{ij}\check h_{ij;}^{\hp{ij;}k}u_ku\\
&=\F\chi^2-\chi^2 \dot\F F^{ij}\{h_{ij;}^{\hp{ij;}k}u_ku+(\chi^{-1}\theta)_ku^kug_{ij}\}\\
&=\F\chi^2-\chi^2\dot\F F^{ij}h_{ij;}^{\hp{ij;}k}u_ku\\
&\q+\dot\F F^{ij}g_{ij}\{\chi_ku^ku\theta-\chi\dot\theta\norm{Du}^2u\}.
\end{aligned}
\end{equation}

Differentiating $\chi$ covariantly with respect to $\xi=(\xi^i)$ we obtain
\begin{equation}
\chi_i=-\chi^2h^k_i\spd{x_k}x,
\end{equation}
\begin{equation}\lae{5.34}
\chi_{ij}=2\chi^{-1}\chi_i\chi_j-\chi^2h^k_{i;j}\spd{x_k}x+h^k_ih_{kj}\chi-\chi^2 h_{ij}.
\end{equation}

Combining \re{5.32} and \re{5.34} the result follows immediately due to the homogeneity of $F$. 
\ep

We want to prove that $h_{ij}$ is uniformly bounded. However, this result can only be achieved in several steps.

We observe that in view of the relation \re{5.11} and the boundedness of $\breve h^i_j$
\begin{equation}
h^i_j (1-\tfrac14u^2)
\end{equation}
is uniformly bounded, or equivalently,
\begin{equation}\lae{5.36}
\abs{h^i_j}e^{-\frac tn}\le \const,
\end{equation}
because of \re{5.4} and \fre{3.17}.

As a first new step we shall improve \re{5.36} slightly:

\bl\lal{5.2}
Define $\lam_\e$ by
\begin{equation}
\lam_\e=\tfrac1n -\e,\qq\e>0,
\end{equation}
where $\e$ is small. Then the principal curvatures $\ka_i$ of the flow hypersurfaces can be estimated from above by
\begin{equation}\lae{5.38}
\ka_i\le c e^{\lam_\e t},
\end{equation}
if $\e>0$ is small
\begin{equation}
0<\e<\e_0.
\end{equation}
\el

\bp
Define $\zeta$ as in \fre{4.32} and let
\begin{equation}
\tilde\zeta =\zeta e^{-\lam_e t},
\end{equation}
then we claim that
\begin{equation}
w=\log \tilde\zeta +\log\chi
\end{equation}
is uniformly bounded from above, if $\e$ is sufficiently small.

Let $0<T<\un$ be large and assume that
\begin{equation}
\sup_{Q_T}w=w(t_0,\xi_0)
\end{equation}
with
\begin{equation}
0<t_0\le T.
\end{equation}

Arguing as in the proof of \frl{4.4}, we may assume that $\ka_n$ is the largest principal curvature and that $w$ is defined by
\begin{equation}
w=\log\tilde h^n_n+\log\chi,
\end{equation}
where
\begin{equation}
\tilde h^n_n=h^n_ne^{-\lam_\e t}.
\end{equation}

We shall suppose that
\begin{equation}
\tilde h^n_n(t_0,\xi_0)>>1.
\end{equation}

Applying the maximum principle we then infer from \re{5.24} and \re{5.26} 
\begin{equation}\lae{5.47}
\begin{aligned}
0&\le \dot\F F^{kl}g_{kl}\{-\tilde\vt v^{-1}e^{\lam_\e t}\tilde h^n_n-v^{-2}\dot{\tilde\vt}+c\dot{\tilde\vt}e^{-\lam_\e t} (\tilde h^n_n)^{-1}+c\tilde\vt\\
&\q+c\dot{\tilde\vt}\norm{Du}^2+c\Ddot{\tilde\vt}\norm{Du}^2e^{-\lam_\e t}(\tilde h^n_n)^{-1}\\
&\q+\tilde\vt u^k(\log\tilde h^n_n)_k+u\theta u^k(\log\chi)_k\},
\end{aligned}
\end{equation}
where we used  the concavity of $F$, the properties of $\F$ and at one point the vanishing of $Dw$ in $(t_0,\xi_0)$.

Our assumption that $\tilde h^n_n$ is very large implies that $t=t_0$ is very large and, hence, powers of $e^t$ will be the dominating terms.

In view of  \fre{3.17} and \re{5.4} we have
\begin{equation}
\tilde\vt\sim ce^{\frac tn}\q\wed\q \dot{\tilde\vt}\sim c e^{\frac{2t}n}\q\wed\q \Ddot{\tilde\vt}\sim ce^{\frac{3t}n},
\end{equation}
while
\begin{equation}
\norm{Du}\le c e^{-\lam_0  t}
\end{equation}
for some $0<\lam_0$, \cf the estimate \fre{3.40}.

The best term inside the braces on the right-hand side of \re{5.47} is
\begin{equation}
-v^{-2}\dot{\tilde\vt}\sim -c e^{\frac {2t}n}
\end{equation}
and the worst is
\begin{equation}
\begin{aligned}
c\Ddot{\tilde\vt}\norm{Du}^2e^{-\lam_e t}(\tilde h^n_n)^{-1}\sim c(\tilde h^n_n)^{-1}e^{\frac tn(2+n\e -2n\lam_0)},
\end{aligned}
\end{equation}
hence, choosing
\begin{equation}\lae{5.52}
\e=2\lam_0
\end{equation}
we obtain an a priori estimate for $\tilde h^n_n$, since the terms in \re{5.47} involving the derivatives of $\log\tilde h^n_n$ and $\log\chi$ vanish, for
\begin{equation}
u\theta =\tilde\vt\q\wed\q Dw=0.
\end{equation}
\ep

As a corollary we deduce:
\bc
The quantity $F=F(\check h^i_j)$ can be estimated from above by
\begin{equation}
F\le nv^{-1}\tilde\vt(1+ce^{-2\lam_0 t})\qq\A\, 0\le t<\un.
\end{equation}
\ec
\bp
This follows at once from \re{5.13}, \re{5.38}, \re{5.52} and the normalization \fre{3.1}.
\ep

We are now able to improve the decay rate of $\abs{Du}$.
\bl\lal{5.4}
$\abs{Du}$ satisfies the estimate
\begin{equation}\lae{5.55}
\abs{Du}\le ce^{-\frac tn}\qq\A\, 0\le t<\un.
\end{equation}
\el
\bp
We look at the scalar flow equation
\begin{equation}\lae{5.56}
\dot u=\pde ut=\frac vF,
\end{equation}
where $F=F(\check h^i_j)$. Let
\begin{equation}
\f=\log u,
\end{equation}
then
\begin{equation}
\begin{aligned}
h^i_j=g^{ik}h_{kj}=v^{-1}u^{-1}\{-(\s^{ik}-v^{-2}\f^i\f^k)\f_{jk}+\de^i_j\},
\end{aligned}
\end{equation}
where all space derivatives are covariant derivatives with respect to $\s_{ij}$, \cf \fre{3.26}, but now we are in euclidean space, i.e., the factor $\vt$ in \re{3.22} is equal to
\begin{equation}
\vt(r)=r.
\end{equation}

Hence, we infer from \re{5.56} 
\begin{equation}
\dot\f=\frac1{F(\tilde h^i_j)},
\end{equation}
where
\begin{equation}
\tilde h^i_j=v^{-2}\{-(\s^{ik}-v^{-2}\f^i\f^k)\f_{jk}+\vt\de^i_j\},
\end{equation}
and $\vt$ is defined by
\begin{equation}
\vt=\vt(r)=r\tilde\vt.
\end{equation}

The term
\begin{equation}
w=\tfrac12\abs{D\f}^2
\end{equation}
then satisfies
\begin{equation}
\begin{aligned}
\dot w&=F^{-2}\{2v^{-2}Fw_i\f^i+v^{-2}F^k_l\tilde g^{lr}w_{kr} -v^{-2}F^k_l\tilde g^{lr}\f_{ik}\f^i_r\\
&\q+v^{-2}F^k_l\tilde g^{lr}_{\hp{lr};i}\f^i\f_{kr}+v^{-2}F^k_l\tilde g^{lr}\f_r\f_k-v^{-2}F^k_l\tilde g^{lr}\s_{kr}\abs{D\f}^2\\
&\q -2v^{-2}F^k_k\dot\vt e^\f w\},
\end{aligned}
\end{equation}
\cf \fre{3.43} observing that now $\vt$ is defined differently and
\begin{equation}
\vt_i\f^i=\dot\vt e^\f\abs{D\f}^2=2\dot\vt e^\f w.
\end{equation}
The metric $\tilde g_{ij}$ is defined by
\begin{equation}
\tilde g_{ij}=\f_i\f_j+\s_{ij},
\end{equation}
and $\tilde g^{ij}$ is its inverse.

Let
\begin{equation}
0<\lam\le\tfrac2n
\end{equation}
be arbitrary and define
\begin{equation}
\tilde w=we^{\lam t}\q\wed\q \tilde\f_i=\f_i e^{\frac{\lam t}2}.
\end{equation}

Now choose $0<T<\un$ and suppose that
\begin{equation}
\sup_{Q_T}w,\qq Q_T=[0,T]\times \Ss[n],
\end{equation}
is attained at $(t_0,x_0)$ with $t_0>0$. Then the maximum principle implies
\begin{equation}\lae{5.70}
\begin{aligned}
0&\le F^{-2}\{F^k_l\tilde g^{lr}\tilde\f_r\tilde\f_k-F^k_l\tilde g^{lr}\s_{kr}\abs{D\tilde\f}^2-2F^k_k\dot\vt e^\f\tilde w+\lam v^2F^2\tilde w\}.
\end{aligned}
\end{equation}

By definition
\begin{equation}
\tilde h^i_j=v^{-1}u^{-1}\check h^i_j,
\end{equation}
and thus
\begin{equation}\lae{5.72}
v^2F(\tilde h^i_j)\le u^2F^2(\check h^i_j)\le u^2n^2v^{-2}\tilde\vt^2(1+ce^{-2\lam_0 t})^2.
\end{equation}

On the other hand, 
\begin{equation}\lae{5.73}
2n\dot\vt u=8n\tilde\vt^2
\end{equation}
hence, we obtain an a priori estimate for $\tilde w$ provided
\begin{equation}
0<\lam<\tfrac2n.
\end{equation}

To derive an a priori estimate in the limit case
\begin{equation}
\lam=\tfrac2n
\end{equation}
we define, with a slight abuse of notation,
\begin{equation}
w=w(t)=\sup_{\Ss[n]}w(t,\cdot)\q\wed\q \tilde w=we^{\lam t}
\end{equation}
with $\lam=\frac{2t}n$; $w$ is Lipschitz continuous and the maximum principle then yields, instead of \re{5.70},
\begin{equation}
\begin{aligned}
\dot{\tilde w}&\le F^{-2}\{-2F^k_k\dot\vt e^\f\tilde w+\lam v^2F^2\tilde w\}
\end{aligned}
\end{equation}
for almost every $t>0$. Because of the relations \re{5.72}, \re{5.73} and the previous estimates for $w$ we then conclude
\begin{equation}
\dot{\tilde w}\le ce^{-\de t}
\end{equation}
for a.e.\ $t>0$ with some $\de>0$, completing the proof of the lemma.
\ep

As a corollary we deduce:
\bc
The principal curvatures $\ka_i$ of the flow hypersurfaces are uniformly bounded from above
\begin{equation}
\ka_i\le c\qq\A\, 0\le t<\un.
\end{equation}
\ec
\bp
Choosing in the proof of \rl{5.2} $\lam_\e=0$ and applying the maximum principle we obtain the inequality  \re{5.47} with  $\lam_\e$ replaced by $0$. In view of the estimate \re{5.55} we then infer an a priori estimate for $\ka_n$.
\ep

An estimate from below for the $\ka_i$ is much more difficult and requires two steps.
\bl
Let $\ka_i$, $1\le i\le n$, be the principal curvatures of the flow hypersurfaces
\begin{equation}
\ka_1\le \cdots\le \ka_n,
\end{equation}
and let
\begin{equation}
\tilde\ka_i(t)=\inf_{\Ss[n]}\ka_i(t,\xi)(2-u(t,\xi)),
\end{equation}
then
\begin{equation}
\liminf_{t\ra\un}\tilde\ka_1(t)=0.
\end{equation}
\el
\bp
We argue by contradiction. Suppose that
\begin{equation}\lae{5.83}
\liminf_{t\ra\un}\tilde\ka_1(t)<0.
\end{equation}
Let $\f$ be defined by
\begin{equation}
\begin{aligned}
\f&=F(2-u)=F(h^i_j+v^{-1}\tilde\vt\de^i_j)(2-u)\\
&=F(h^i_jv\tilde\vt^{-1}+\de^i_j)v^{-1}\tilde\vt (2-u),
\end{aligned}
\end{equation}
then
\begin{equation}
\lim_{t\ra\un}\,\abs{v^{-1}\tilde\vt (2-u)}_{0,\Ss[n]}=1 
\end{equation}
and \re{5.83} is equivalent to
\begin{equation}
\liminf_{t\ra\un}\inf_{\Ss[n]}F(2-u)=F(1+\tilde\ka_1,\ldots,1+\tilde\ka_n)<F(1,\ldots,1)=n,
\end{equation}
since the non-negative $\ka_i$ are uniformly bounded and $F$ is strictly monotone. Thus, \re{5.83} implies
\begin{equation}
\liminf_{t\ra\un}\inf_{\Ss[n]}w<\log n,
\end{equation}
where
\begin{equation}
\begin{aligned}
w&=\log\f-\log\chi-\log 2\\
&=-\log(-\F)+\log(2-u)-\log\chi-\log 2,
\end{aligned}
\end{equation}
since
\begin{equation}\lae{5.89}
\lim_{t\ra\un}\,\abs{-\log\chi-\log 2}_{0,\Ss[n]}=0.
\end{equation}

Let $\e>0$ be so small such that
\begin{equation}\lae{5.90}
\liminf_{t\ra\un}\inf_{\Ss[n]}w<(1-2\e)\log n
\end{equation}
and let $\tau$ be so large such that
\begin{equation}
t\ge \tau\q\wed\q \inf_{\Ss[n]}w(t,\xi)=w(t,\xi_0)<(1-\e)\log n
\end{equation}
implies
\begin{equation}
\tilde\ka_1(t)<-\e_0\q\wed\q\ka_1\tilde\vt^{-1}(t,\xi_0)<-\e_0
\end{equation}
for a fixed $0<\e_0=\e_0(\e)$. The existence of $\tau$ follows from the relations \re{5.89} and \re{5.90}.

Define the set
\begin{equation}
\Lam=\set{t}{t\ge\tau\q\wed\q \inf_{\Ss[n]}w(t,\cdot)<(1-\e)\log n},
\end{equation}
then $\Lam\ne\eS$, since it contains a sequence $t_k\ra\un$.

We shall now prove
\begin{equation}\lae{5.94}
\Lam=[\tau,\un),
\end{equation}
and
\begin{equation}
\tilde w(t)=\inf_{\Ss[n]}w(t,\cdot)
\end{equation}
is (weakly) monotone increasing in $[\tau,\un)$, i.e.,
\begin{equation}\lae{5.96} 
\tilde w(t_1)\le \tilde w(t_2)\qq\A\,\tau\le t_1<t_2<\un.
\end{equation}

Let $T$, $\tau<T<\un$, be arbitrary but so large such that 
\begin{equation}
\Lam\ii [\tau,T]\ne\eS,
\end{equation}
and suppose that
\begin{equation}
\inf\set{w(t,\xi)}{\tau\le t\le T,\; \xi\in \Ss[n]}=w(t_0,\xi_0)
\end{equation}
with $t_0>\tau$.

$w$ satisfies the evolution equation
\begin{equation}\lae{5.99}
\begin{aligned}
&\msp[140]\dot w-\dot\F F^{ij}w_{ij}=\\
&\;\dot\F F^{ij}g_{ij}r_{\al\bet}\nu^\al\nu^\bet\tilde\vt+\dot\F F^{ij}g_{ij}\dot{\tilde\vt}v^{-2}
 -\dot\F F^{ij}g_{ij}(\log(-\F))_ku^k\tilde\vt \\
&-\dot\F F^{ij}(\log(-\F))_i(\log(-\F))_j+\dot\F F^{ij}(\log(2-u))_i(\log(2-u))_j\\
&+\dot\F F^{ij}(\log\chi)_i(\log\chi)_j\\
&+\dot\F F^{ij}g_{ij}\tilde\vt v^{-2}(2-u)^{-1}+\dot\F F^{ij}\bar h_{ij}(2-u)^{-1}-2v^{-1}\f^{-1}\\
&+\dot\F F^{ij}g_{ij}\{\tilde\vt\chi -(\log\chi)_ku^ku\theta+\dot\theta \norm{Du}^2u\}.
\end{aligned}
\end{equation}

At the point $(t_0,\xi_0)$ $Dw=0$, hence the terms in line two and three on the right-hand side of \re{5.99} add up to
\begin{equation}
2\dot\F F^{ij}(\log\chi)_i(\log(2-u))_j,
\end{equation}
which in turn is equal to
\begin{equation}
2\dot\F F^{ij}h_{ik}u^ku_j(2-u)^{-1}\chi^{-1}\ge -c\dot\F,
\end{equation}
due to the estimates \re{5.55}, \re{5.36}, \re{5.4} and \fre{3.17}.

Analogously, we conclude
\begin{equation}
\begin{aligned}
&\dot\F F^{ij}g_{ij}\{-(\log(-\F))_ku^ku\theta-(\log\chi)_ku^k\tilde\vt\}=\\
&\qq \dot\F F^{ij}g_{ij}(-\log(2-u))_ku^k)\tilde\vt\ge0,
\end{aligned}
\end{equation}
where we also used 
\begin{equation}
u\theta=\tilde\vt.
\end{equation}

Hence, applying the maximum principle we infer from \re{5.99} 
\begin{equation}\lae{5.104}
\begin{aligned}
0&\ge \dot\F F^{ij}g_{ij}\dot{\tilde\vt}v^{-2}+\dot\F F^{ij}g_{ij}\tilde\vt v^{-2}(2-u)^{-1}-2v^{-1}\f^{-1}-c\dot\F\\
&\ge 2\de-ce^{-\frac{2t}n}\ge \de,
\end{aligned}
\end{equation}
if $\tau$ is large enough, with some uniform $\de=\de(\e_0)>0$; a contradiction.

Thus, we have proved that $t_0=\tau$ and therefore 
\begin{equation}
\tilde w(\tau)\le \tilde w(t)\qq\A\, \tau\le t\le T.
\end{equation}

Since we can replace $\tau$ by any $t_1\in[\tau,T)$ we conclude
\begin{equation}
\tilde w(t_1)\le \tilde w(t_2)\qq\A\, \tau\le t_1\le t_2\le T,
\end{equation}
and we have proved \re{5.94} as well as \re {5.96}, since $\tau<T<\un$ is arbitrary.

However, the arguments we used to derive the contradiction in inequality \re{5.104} yield
\begin{equation}
\dot w(t,\xi_t)\ge \de >0\qq\A\, \tau\le t<\un,
\end{equation}
where
\begin{equation}
\inf\set{w(t,\xi)} {\xi\in \Ss[n]}=w(t,\xi_t),
\end{equation}
in view of \re{5.94} and the definition of $\Lam$. The left-hand side of the preceding equation is the definition of $\tilde w(t)$, which is Lipschitz continuous and satisfies for a.e.\ $t>\tau$
\begin{equation}
\dot{\tilde w}(t)=\dot w(t,\xi_t).
\end{equation}
Hence, we deduce
\begin{equation}
\dot{\tilde w}(t)\ge \de
\end{equation}
for a.e.\  $t>\tau$, which is a contradiction, since $\tilde w$ is uniformly bounded, completing the proof of the lemma.
\ep

Now, we can prove that the principal curvatures are uniformly bounded.

\bl
The principal curvatures $\ka_i$, $1\le i\le n$, are uniformly bounded during the evolution
\begin{equation}\lae{5.111}
\abs{\ka_i}\le c.
\end{equation}
\el
\bp
We shall estimate
\begin{equation}
\f=\tfrac12\abs A^2=\tfrac12 h_{ij}h^{ij},
\end{equation}
which satisfies the evolution equation
\begin{equation}\lae{5.113}
\dot\f-\dot\F F^{ij}\f_{ij}=-\dot\F F^{kl}h_{ij;k}h^{ij}_{\hp{ij};l}+\{\dot h^j_i-\dot\F F^{kl}h^j_{i;kl}\}h^i_j.
\end{equation}
Looking at \re{5.24} and observing that, in view of the previous estimates,
\begin{equation}
\lim_{t\ra\un}\abs Ae^{-\frac tn}=0
\end{equation}
uniformly in $\xi\in \Ss[n]$, and, because of the homogeneity of $F$,
\begin{equation}
F^{-1}+\abs{F^{kl,rs}}\le c e^{-\frac tn},
\end{equation}
we deduce that the terms
\begin{equation}
-\dot\F F^{kl}h_{ij;k}h^{ij}_{\hp{ij};l}-\dot\F F^{kl}g_{kl}v^{-2}\dot{\tilde\vt}2\f,
\end{equation}
which are either explicitly or implicitly contained in the right-hand side of \re{5.113}, are dominating; they can absorb any bad term such that an application of the maximum principle gives an a priori estimate for $\f$.
\ep

As a corollary we obtain:
\bt
The flow hypersurfaces in hyperbolic space become strongly convex exponentially fast and also more and more totally umbilic. In fact there holds
\begin{equation}\lae{5.117}
\abs{\breve h^i_j-\de^i_j}\le c e^{-\frac tn}.
\end{equation}
\et
\bp
We infer from \fre{5.11}
\begin{equation}
\begin{aligned}
\breve h^i_j-\de^i_j&=\breve h^i_j-\frac u{2v}\de^i_j+(\frac u{2v}-1)\de^i_j\\
&= h^i_j\frac u{2\tilde\vt}+(\frac u{2v}-1)\de^i_j,
\end{aligned}
\end{equation}
from which the estimate \re{5.117} immediately follows, in view of \re{5.111}, \re{5.55}, \re{5.4}, and \fre{3.17}. 
\ep

\section{Higher order estimate}
Assuming the curvature function $F$ to be smooth, we want to prove higher order estimates for $h_{ij}$, or equivalently, for $u$. Since we already know that $h_{ij}$ is uniformly bounded,
\begin{equation}
\dot g_{ij}=-2\F h_{ij}
\end{equation}
as well as the Riemannian curvature tensor of the induced metric then are also uniformly bounded. 

Let $A$ represent the second fundamental form, where we omit the tensor indices, then we want to prove
\begin{equation}\lae{6.2}
\norm{D^mA}\le c_m e^{-\frac tn}\qq\A\, m\ge 1.
\end{equation}

This estimate will immediately imply a corresponding estimate
\begin{equation}
\norm{D^mu}\le c_m e^{-\frac tn}\qq\A\, m\ge 1,
\end{equation}
in view of the relation
\begin{equation}\lae{6.4}
h_{ij}v^{-1}=-u_{ij}+\bar h_{ij}
\end{equation}
and the estimate \fre{5.55}.

To obtain an estimate for $D^2u$ we have to apply the following interpolation lemma:

\bl
Let $M=M^n$ be a compact Riemannian manifold of class  $C^m$, $m\ge 2$, and $u\in C^2(M)$, then
\begin{equation}
\norm{Du}_{0,M}\le c \msp\abs u_{2,M}^\frac12\msp \abs u_{0,M}^\frac12,
\end{equation}
where $c=c(M)$ and
\begin{equation}
\abs u_{2,M}=\abs u_{0,M}+\norm{Du}_{0,M}+\norm{D^2u}_{0,M}
\end{equation}
and the norms on the right-hand side are supremum norms.
\el
\bp
Using a partition of unity we may assume that $u$ has support in a coordinate chart and hence we may assume that
\begin{equation}
u\in C^2_c(\R[n]),
\end{equation}
where $\R[n]$ is equipped with the Euclidean metric. Moreover, we may assume $n=1$.

Let $x\in \R[]$, $\e>0$ be arbitrary, and choose $x_1, x_2\in\R[]$ such that
\begin{equation}
x_2-x_1=\e\q\wed\q x\in (x_1,x_2).
\end{equation}
Then we deduce
\begin{equation}
u(x_2)-u(x_1)=Du(\xi)(x_2-x_1),\q \xi\in (x_1,x_2),
\end{equation}
\begin{equation}
Du(x)=Du(\xi)+\int_\xi^xD^2u,
\end{equation}
and hence,
\begin{equation}
\abs{Du(x)}\le 2\e^{-1}\abs u_0+\e\abs{D^2u}_0\equiv \f(\e),
\end{equation}
where we assume without loss of generality that
\begin{equation}
\abs{D^2u}_0>0,
\end{equation}
otherwise, we replace $\abs{D^2u}_0$ by $\abs{D^2u}_0+\de$, $\de>0$.

Minimizing $\f$ by solving
\begin{equation}
\dot\f(\e)=-2\e^{-2}\abs u_0+\abs{D^2u}_0=0,
\end{equation}
we conclude
\begin{equation}
\e=\sqrt 2\msp[1]\abs u_0^\frac12\msp\abs{D^2u}_0^{-\frac12},
\end{equation}
and thus,
\begin{equation}
\abs{Du}\le 2 \sqrt 2\msp \abs u_0^\frac12\msp \abs{D^2u}_0^\frac12.
\end{equation}
\ep

\bc\lac{6.2}
Let $M=M^n$ be a compact Riemannian manifold of class $C^m$, $m\ge 2$, and $u\in C^m(M)$, then
\begin{equation}
\norm{D^{m-1}u}_{0,M}\le c\msp \abs u _{m-2,M}^\frac12\msp \abs u_{m,M}^\frac12,
\end{equation}
where $c=c(m,M)$.
\ec

We shall apply the corollary to the function $(u-2)$ using either $M=\Ss[n]$ or $M=\graph u$.

The starting point for deriving the estimate \re{6.2} is equation \fre{5.24} which will be differentiated covariantly. However, we first have to derive some preparatory lemmata. 

\bl
Let $\f$ be defined by
\begin{equation}
\f=(u-2)^{-1}
\end{equation}
and assume 
\begin{equation}\lae{6.18}
m\ge 2\q\wed\q \norm{D^kA}\le c_\lam e^{-\lam t} \q\A\, 1\le k\le m-1,
\end{equation}
and for all $0<\lam<\frac1n$, then
\begin{equation}\lae{6.19}
\norm{D^{m+1}\f}\le c_\e \msp e^{(\frac1n+\e)t}\qq \A\,  0<\e.
\end{equation}
The estimate
\begin{equation}\lae{6.20}
\norm{D\f}\le c\msp e^{\frac tn}
\end{equation}
has already been proved.
\el
\bp
Set 
\begin{equation}
\tilde u=(u-2)\msp e^\frac tn,
\end{equation}
then $\tilde u$ satisfies
\begin{equation}
-c_1\le \tilde u\le -c_2\qq\A\, 0\le t<\un,
\end{equation}
and 
\begin{equation}
\f=\tilde u^{-1}e^\frac tn.
\end{equation}
Let $\al\in\N^n$ be a multi-index of order $m+1$, $m\ge 2$,  then $D^\al\f$ can be written as
\begin{equation}\lae{6.24}
D^\al\f=\sum_{\abs{\bet_1}+\cdots+\abs{\bet_{m+1}}=m+1}c_{\bet_1,\ldots,\bet_{m+1}}D^{\bet_1}\tilde u\cdots D^{\bet_{m+1}}\tilde u\msp e^\frac tn,
\end{equation}
where the coefficients $c_{\bet_1,\ldots,\bet_{m+1}}$ depend smoothly on $\tilde u$, and, if we allow some of the coefficients to vanish, the sum is taken over all multi-indices $\bet_i$, $1\le i\le m+1$, satisfying
\begin{equation}
\sum_{i=1}^{m+1}\abs{\bet_i}=m+1.
\end{equation}

The estimate \re{6.20} is trivial in view of \fre{5.55}.
\ep

\bl\lal{6.4}
Let $f=f(u,Du,\tilde u, D\tilde u)$ be any smooth function and assume that the conditions \re{6.18} are valid, then, for any $\e>0$, there holds
\begin{equation}\lae{6.26}
\norm{D^m(f\tilde\vt)}\le c_\e\msp e^{(\frac1n+\e)t},
\end{equation}
\begin{equation}
\norm{D^m(f\dot{\tilde\vt}u_i)}\le c_\e\msp e^{(\frac1n+\e)t},\q\A\, 1\le i\le n,
\end{equation}
\begin{equation}
\norm{D^m(f\Ddot{\tilde\vt}u_iu_j)}\le c_\e\msp e^{(\frac1n+\e)t},\q\A\, 1\le i,j\le n.
\end{equation}
\el
\bp
Let us only prove \re{6.26}, since we can write
\begin{equation}
\begin{aligned}
\dot{\tilde\vt}u_i&=\dot{\tilde\vt}(u-2)(u-2)^{-1} u_i\\
&=\dot{\tilde\vt}\tilde u^{-1}\tilde u_i\\
&\equiv f\tilde\vt
\end{aligned}
\end{equation}
with some smooth function $f=f(u,Du,\tilde u,D\tilde u)$, and similarly
\begin{equation}
\Ddot{\tilde\vt}u_iu_j\equiv f\tilde\vt.
\end{equation}

\cvm
\cq{\re{6.26}}\q   Define $\tilde\theta$ by
\begin{equation}
\tilde\theta=-\tilde\vt (u-2),
\end{equation}
then $\tilde\theta$ is smooth and
\begin{equation}
f\tilde\vt=-f\tilde\theta\f.
\end{equation}
The estimate then follows by applying the general Leibniz rule and \re{6.19} observing that $\e>0$ is assumed to be arbitrary.
\ep

Let $\check A$ be a symbol for $\check h^i_j$, then
\begin{equation}
\check A=A+v^{-1}\tilde\vt\de^i_j
\end{equation}
and we deduce from \re{6.26}, if the assumptions \re{6.18} are satisfied,
\begin{equation}\lae{6.34}
\norm{D^m\check A}\le \norm{D^mA}+c_\e\msp e^{(\frac1n+\e)t}
\end{equation}
for any $\e>0$, where $c_\e$ also depends on $m$. Here, we also used the relation
\begin{equation}
v^{-2}=1-\norm{Du}^2.
\end{equation}
We also note that in case $m=1$ the relation \re{6.34} is valid for $\e=0$.

Inside the braces of the right-hand side of equation \fre{5.24} there is the crucial term
\begin{equation}\lae{6.36}
-v^{-2}\dot{\tilde\vt}h^j_i+v^{-1}\dot{\tilde\vt}\bar h_{ik}g^{kj},
\end{equation}
which is equal to
\begin{equation}
v^{-1}\dot{\tilde\vt}\{-v^{-1}h^j_i+\bar h_{ik}g^{kj}\}=v^{-1}\dot{\tilde\vt}u^j_i,
\end{equation}
in view of \re{6.4}. 

Differentiating \re{6.36} covariantly with respect to a multi-index $\al$, $\abs\al=m$, $m\ge 1$, we therefore obtain
\begin{equation}\lae{6.38}
\begin{aligned}
&\sum_{\bet\le\al}\binom\al\bet D^{\al-\bet}(v^{-1}\dot{\tilde\vt})D^\bet\{-v^{-1}h^j_i+\bar h_{ik}g^{kj}\}\\
&=v^{-1}\dot{\tilde\vt}D^\al\{-v^{-1}h^j_i+\bar h_{ik}g^{kj}\}+\sum_{\bet<\al}\binom\al\bet D^{\al-\bet}(v^{-1}\dot{\tilde\vt})D^\bet u^j_i.
\end{aligned}
\end{equation}

Furthermore, there holds
\begin{equation}
v^{-1}\dot{\tilde\vt}=f(u,Du,\tilde u)\tilde\vt e^\frac tn,
\end{equation}
and hence,
\begin{equation}
D^{\al-\bet}(v^{-1}\dot{\tilde\vt})D^\bet u^j_i=D^{\al-\bet}(f\tilde\vt)D^\bet\tilde u^j_i
\end{equation}
and we conclude
\begin{equation}\lae{6.41}
\begin{aligned}
\norm{\sum_{\bet<\al}\binom\al\bet D^{\al-\bet}(v^{-1}\dot{\tilde\vt})D^\bet u^j_i}\le c_\e\msp e^{(\frac1n+\e)t}\q\A\, \e>0,
\end{aligned}
\end{equation}
provided \re{6.18} is valid, in view of \re{6.26}, \re{6.24} and \re{6.19}, since
\begin{equation}
\bet<\al\q\im \q\abs\bet +2\le m+1.
\end{equation}

In case $m=1$ we have
\begin{equation}\lae{6.43}
\norm{D(v^{-1}\dot{\tilde\vt})u^j_i}\le c\msp e^{\frac 2nt}\msp \norm{D^2 u}.
\end{equation}

Hence we obtain:
\bl\lal{6.5}
Let $\al$ be a multi-index of order $m\ge 2$ and suppose that \re{6.18} is valid, then
\begin{equation}\lae{6.44}
\begin{aligned}
D^\al\{-v^{-2}\dot{\tilde\vt}h^j_i+v^{-1}\dot{\tilde\vt}\bar h_{ik}g^{kj}\}&=v^{-1}\dot{\tilde\vt}D^\al\{-v^{-1}h^j_i+\bar h_{ik}g^{kj}\}+\mc O^1_\e\\
&=-v^{-2}\dot{\tilde\vt}D^\al h^j_i+\mc O^1_\e,
\end{aligned}
\end{equation}
where $\mc O^1_\e$ represents a tensor that can be estimated like
\begin{equation}
\norm{\mc O^1_\e}\le c_\e\msp e^{(\frac1n+\e)t}\q\A\, \e>0.
\end{equation}
In case $m=1$ we have
\begin{equation}\lae{6.46}
\begin{aligned}
D\{-v^{-2}\dot{\tilde\vt}h^j_i+v^{-1}\dot{\tilde\vt}\bar h_{ik}g^{kj}\}&=v^{-1}\dot{\tilde\vt}D\{-v^{-1}h^j_i+\bar h_{ik}g^{kj}\}+\mc O^2_0u^j_i\\
&=-v^{-2}\dot{\tilde\vt}Dh^j_i+\mc O^1_0+\mc O^2_0u^j_i,
\end{aligned}
\end{equation}
where $\mc O^1_0$ \resp $\mc O^2_0$ represent tensors that can be estimated like 
\begin{equation}
\norm{\mc O^1_0}\le c\msp e^\frac tn
\end{equation}
\resp
\begin{equation}
\norm{\mc O^2_0}\le c\msp e^{\frac2n t}.
\end{equation}
\el
\bp
Observing that
\begin{equation}
\bar h_{ik}=u^{-1}\bar g_{ik}=u^{-1}g_{ik}-u^{-1}u_iu_k
\end{equation}
the relation \re{6.44} follows from \re{6.38}, \re{6.41} and \re{6.18}, while \re{6.46} can be deduced from \re{6.43}.
\ep

\bd\lad{6.6}
Let $k\in\Z$, then the symbol $\mc O^k_\e$ represents any tensor that can be estimated by
\begin{equation}
\norm{\mc O^k_\e}\le c_\e\msp e^{(\frac kn+\e)t}\qq\A\,\e>0,
\end{equation}
and the symbol $\mc O^k_0$ represents any tensor that can be estimated by
\begin{equation}
\norm{\mc O^k_0}\le c\msp e^{\frac knt}.
\end{equation}
Thus $\mc O^0_0$ represents a uniformly bounded tensor.
\ed
We also denote by $\mc D^mF$ the derivatives of order $m$ of $F$ with respect to the argument $\check h^i_j$, and when $S$, $T$ are arbitrary tensors then $S\star T$ will symbolize any linear combination of tensors formed by contracting over $S$ and $T$. The result can be a tensor or a function. Note that we do not distinguish between $S\star T$ and $cS\star T$, $c$ a constant.

From \re{6.34}, the homogeneity of $F$ and the definition of $\F$ we then deduce
\bl\lal{6.7}
Let $m=1$ or assume that \re{6.18} is valid, then we have
\begin{equation}
D^mF=\mc DF\star D^mA+\mc DF\star \mc O^1_\e,
\end{equation}
\begin{equation}
D^m\F^{(k)}=\F^{(k+1)}D^mF+\F^{(k+1)}\star \mc O^1_\e,
\end{equation}
and similarly
\begin{equation}
D^m\mc D^kF=\mc D^{k+1}F\star D^mF+\mc D^{k+1}F\star\mc O^1_\e,
\end{equation}
where $\F^{(k)}$ is the $k$-th derivative of $\F$. In case $m=1$ $\mc O^1_\e$ can be replaced by $\mc O^1_0$.
\el
We are now ready to differentiate \fre{5.24} covariantly. 
\bl
The tensor $DA$ satisfies the evolution equation 
\begin{equation}\lae{6.55}
\begin{aligned}
&\qq\qq\qq\frac D{dt}(DA)-\dot\F F^{kl}(DA)_{;kl}=\\
&\msp[6] \Ddot\F \,\mc O^0_0\star (DA+\mc O^1_\e)\star\mc D^2A+\dot\F\,\mc D^2F\star(DA+\mc O^1_\e)\star D^2A\\
&+\dot\F\,\mc O^0_0\star DA+\Ddot\F\, \mc O^0_0\star (DA+\mc O^1_0)+\dot\F\,\mc O^0_0\star(DA+\mc O^1_0)\\
&+\F\,\mc O^0_0\star DA+\dddot{\msp[-2]\F}\,\mc O^0_0\star(DA+\mc O^1_0)\star (DA+\mc O^1_0)\star (DA+\mc O^1_0)\\
&+\ddot\F\,\mc O^0_0\star(D^2A+\mc O^1_0+\mc O^2_0\star D^2u)\star(DA+\mc O^1_0)\\
&+\ddot\F\,\mc D^2F\star\mc DF\star (DA+\mc O^1_0)\star (DA+\mc O^1_0)\star (DA+\mc O^1_0)\\
&+\dot\F\, \mc D^3F\star (DA+\mc O^1_0)\star (DA+\mc O^1_0)\star (DA+\mc O^1_0)\\
&+\dot\F\, \mc D^2F\star (D^2A+\mc O^1_0+\mc O^2_0\star D^2u)\star (DA+\mc O^1_0)\\
&+\ddot\F\, \mc O^0_0\star (DA+\mc O^1_0)\star (\dot{\tilde\vt}D^2u\star \mc O^0_0+\tilde\vt\,\mc O^0_0+DA\star\mc O^0_0)\\
&+\dot\F\,\mc O^0_0\star(\mc O^1_0+\mc O^1_0\star DA+\mc O^0_0\star D^2A+\mc O^2_0\star D^2u)\\
&+\dot\F F^{kl}g_{kl}(-v^{-2}\dot{\tilde\vt}DA).
\end{aligned}
\end{equation}
\el
\bp
Differentiate \fre{5.24} covariantly with respect to a spatial variable and apply \frl{5.4}, \rl{6.5}, \rd{6.6},  and \rl{6.7}.
\ep
An almost identical proof---where we also have to rely on \rl{6.4}---yields the evolution equation for higher derivatives of $A$.
\bl
Let $m\ge 2$ and assume that assumptions \re{6.18} are valid, then the tensor $D^mA$, where $D^mA$ represents any covariant derivative $D^\al A$, $\abs\al=m$, satisfies the evolution equation
\begin{equation}\lae{6.56}
\begin{aligned}
&\qq\qq\qq\q\frac D{dt}(D^mA)-\dot\F F^{kl}(D^mA)_{;kl}=\\
&\msp[6]\ddot\F\,\mc DF\star (D^mA+\mc O^1_\e)\star D^2A+\dot\F\,\mc D^2F\star (D^mA+\mc O^1_\e)\star D^2A\\
&+\F\,\mc O_0^0\star D^mA+\dot\F\,\mc O^1_\e\star D^{m+1}A+\dot\F\,\mc D^2F\star\mc O^1_\e\star (D^{m+1}A+\mc O^1_\e)\\
&+\ddot\F\,\mc D^2F\star(D^2F+\mc O^1_\e)\star\mc DF\star(D^mA+\mc O^1_\e)\\
&+\dot\F\, \mc D^2F\star (D^2A+\mc O^1_\e)\star (D^mA+\mc O^1_\e)\star \mc DF\\
&+\ddot\F\,(\mc DF\star D^mA+\mc O^1_\e)\star \mc O^0_0+\dot\F\,\mc D^2F\star  (D^mA+\mc O^1_\e)\star \mc O^0_0\\
&+\dot\F\,(D^mA+\mc O^1_\e)\star \mc O^0_0+\dddot{\msp[-2]\F}(D^mA+\mc O^1_\e)\star D\check A\star D\check A\star\mc O^0_0\\
&+\ddot\F\,\mc D^2F\star D^mA\star D\check A+\ddot\F\,\mc D^2F\star \mc O^1_\e\star D\check A\\
&+\ddot\F\,\mc DF\star(D^mA+\mc O^1_\e)\star\mc D^2F\star D\check A\star D\check A\star\mc O^0_0\\
&+\dot\F\,\mc D^3F\star (D^mA+\mc O^1_\e)\star D\check A\star D\check A\\
&+\dot\F\, \mc D^2F\star (D^{m+1}A+\mc O^1_\e)\star (DA+\mc O^1_\e)\\
&+\ddot\F\, (D^mA+\mc O^1_\e)\star (\dot{\tilde\vt}D^2u\star \mc O^0_0+\tilde\vt \,\mc O^0_0)\\
&+\dot\F\,\mc D^2F\star (D^mA+\mc O^1_\e) \star(\dot{\tilde\vt}D^2u+\tilde\vt\,\mc O^0_0)\star\mc O^0_0\\
&+\dot\F\, \mc DF\star(\mc O^1_\e+\tilde\vt D^mA\star\mc O^0_0+\mc O^0_0\star D^{m+1}A)\\
&-\dot\F F^{kl}g_{kl}v^{-2}\dot{\tilde\vt}D^mA
\end{aligned}
\end{equation}
\el
We are now going to prove uniform bounds for
\begin{equation}
\tfrac12\norm{D^m\tilde A}^2=\tfrac12\sum_{\abs\al=m}\norm{D^\al A}^2e^{2\lam t}
\end{equation}
for all $m\ge 1$ and
\begin{equation}
0\le \lam <\tfrac1n.
\end{equation}

First, we observe that
\begin{equation}\lae{6.59}
\begin{aligned}
&\qq\qq\frac D{dt}(\tfrac12\norm{D^m\tilde A}^2)-\dot\F F^{kl} (\tfrac12\norm{D^m\tilde A}^2)_{;kl}=\\
&\Big\{\frac D{dt}(D^mA)-\dot\F F^{kl}(D^mA)_{;kl}\Big\}e^{\lam t}D^m\tilde A-\dot\F F^{kl}(D^m\tilde A)_{;k}(D^m\tilde A)_{;l}\\
&\qq\qq\qq\qq+\lam \norm{D^m\tilde A}^2
\end{aligned}
\end{equation}
\bl
The quantities $\frac12 \norm{D^m\tilde A}^2$ are uniformly bounded during the evolution for any $m\ge 1$ and $0\le \lam <\frac1n$,
\el
\bp
We prove the lemma recursively by estimating
\begin{equation}
\f=\log(\tfrac12\norm{D^m\tilde A}^2)+\mu\tfrac12\norm{D^{m-1}A}^2,
\end{equation}
where
\begin{equation}
0<\mu=\mu(m)<<1,
\end{equation}
\cf the proof of \cite[Lemma 7.6.3]{cg:cp}.

We shall only treat the case $m=1$, since the proof for $m\ge 2$ is almost identical by considering the evolution equation \re{6.56} instead of \re{6.55}.

Thus, let
\begin{equation}
\f=\log(\tfrac12\norm{D\tilde A}^2)+\mu\tfrac12\norm{A}^2.
\end{equation}
Fix $0<T<\un$, $T$ very large, and suppose that
\begin{equation}
\sup_{[0,T]}\sup_{M(t)}\f=\f(t_0,\xi_0)
\end{equation}
is large, and hence, $0<t_0\le T$, is sufficiently large, such that the previous decay estimates for $\norm{Du}$, etc.\  can be employed.

Applying the maximum principle we deduce from \re{6.55}, \re{6.59} and the evolution equation for $\frac12\norm{A}^2$, see \fre{5.113},
\begin{equation}\lae{6.64}
\begin{aligned}
0&\le \{-\dot\F F^{kl}g_{kl}v^{-2}\dot{\tilde\vt}+\lam\}-2\dot\F F^{kl}(D\tilde A)_{;k}(D\tilde A)_{;l}\\
&\q+\dot\F F^{kl}\log(\tfrac12\norm{\tilde A}^2)_k\log(\tfrac12\norm{D\tilde A}^2)_l\\
&-\mu\dot\F F^{kl}h_{ij;k}h^{ij}_{\hp{ij};l}+\text{rest}.
\end{aligned}
\end{equation}

When $t_0$ is large then terms in the braces can be estimated from above by
\begin{equation}
-2\de,
\end{equation}
where 
\begin{equation}
\de=\de(\lam)\approx \tfrac12 (\tfrac1n-\lam).
\end{equation}

To estimate
\begin{equation}
\dot\F F^{kl}\log(\tfrac12\norm{\tilde A}^2)_k\log(\tfrac12\norm{D\tilde A}^2)_l
\end{equation}
we use $D\f=0$ and conclude that this term can be estimated from above by
\begin{equation}
\begin{aligned}
\mu^2\dot\F F^{kl}h_{ij;k}h^{ij}_{\hp{ij};l}\norm{A}^2\le \tfrac\mu 2\dot\F F^{kl}h_{ij;k}h^{ij}_{\hp{ij};l},
\end{aligned}
\end{equation}
if $0<\mu$ is small.

Most terms in the \cq{rest} can be easily absorbed; a few are a bit more delicate. These can be estimated from above by
\begin{equation}
c\msp\frac{\norm{D^2u}e^{\lam t}\norm{D\tilde A}}{\norm{D\tilde A}^2}<c\msp\norm{D\tilde A}^{-\frac12},
\end{equation}
where the last inequality is due to the interpolation lemma, \cf \rc{6.2}, applied to $(u-2)$. Here, we also used the assumption that $\norm{D\tilde A}\ge 1$.

A thorough inspection of the right-hand side of \re{6.64} then yields
\begin{equation}
0\le -\de+c\msp e^{-(\frac1n-\lam)t_0}+c\msp\norm{D\tilde A}^{-\frac12}
\end{equation}
and hence an a priori estimate for $\norm{D\tilde A}$, if $t_0$ is large.
\ep
It remains to prove the optimal decay \re{6.2} and the convergence to a constant. The optimal decay will be achieved by deriving the equivalent estimate
\begin{equation}\lae{6.71}
\norm{D^mu}\le c_me^{-\frac tn}\qq\A\,m\ge 1.
\end{equation}
\bt
Let $M(t)=\graph u(t)$ be the leaves of the inverse curvature flow, where $F$ and the initial hypersurface are smooth, then the estimate \re{6.71} is valid and the function
\begin{equation}\lae{6.72}
(u-2)e^{\frac tn}
\end{equation}
converges in $C^\un(\Ss[n])$ to a strictly negative function.
\et
\bp
It suffices to prove \re{6.71}, in view of the relations \fre{3.17} and \fre{5.4}, and to show that the limit exists.

\cvm
(i) Our starting point is equation \fre{5.21} satisfied by $u$ as well as by $(u-2)$. 

Let $\f$, $\tilde \f$, $\tilde F$, and $\tilde \F$ be defined by
\begin{equation}
\f=(2-u)^{-1}\q\wed\q \tilde\f=\f e^{-\frac tn},
\end{equation}
\begin{equation}
\tilde F=F(\check h^k_l(2-u)),
\end{equation}
and 
\begin{equation}
\tilde \F=\F(\tilde F),
\end{equation}
then we deduce from \re{5.21}
\begin{equation}\lae{6.76}
\begin{aligned}
\dot{\tilde\f}-\dot{\tilde\F}\tilde\f^{-2}e^{-\frac2n t}\tilde F^{ij}\tilde\f_{ij}&=-2\dot{\tilde\F}\tilde\f^{-2}e^{-\frac2nt}\tilde F^{ij}\tilde\f_i\tilde \f_j\tilde \f^{-1}+2v^{-1}\tilde F^{-1}\tilde\f\\
&\hp{=}\msp[6] -\tilde F^{-2} F^{ij}g_{ij}\tilde\theta \tilde\f v^{-2}-\tilde F^{-2} F^{ij}\bar h_{ij}e^{-\frac tn}-\tfrac1n \tilde \f,
\end{aligned}
\end{equation}
where
\begin{equation}\lae{6.77}
\tilde\theta=\tilde\vt(2-u).
\end{equation}
$\tilde\theta$ depends smoothly on $u$ and is strictly positive. The derivatives of arbitrary order of $\tilde\theta$, $\tilde F$, $\tilde\F$, $F^{ij}$, $\tilde F^{ij}$, $v$, and $\bar h_{ij}$ are uniformly bounded and decay expo\-nentially fast, if $t$ goes to infinity, while the $C^m$-norms of $\tilde\f$ can be estimated by
\begin{equation}\lae{6.78}
\norm{D^m\tilde\f}\le c_{m,\e}\msp e^{\e t}\qq\A\, \e>0,
\end{equation}
in view of our previous estimates.

Differentiating then \re{6.76} covariantly we obtain the following differential inequality
for
\begin{equation}
w=\tfrac12\norm{D^m\tilde\f}^2
\end{equation}
\begin{equation}\lae{6.80}
\begin{aligned}
\dot w -\dot{\tilde\F}\tilde\f^{-2}e^{-\frac2n t}\tilde F^{ij}w_{ij}\le \mc O^0_{-2\de}+2\{2v^{-1}\tilde F^{-1}-\tilde F^{-2}F^{ij}g_{ij}\tilde\theta v^{-2}-\tfrac1n\}w,
\end{aligned}
\end{equation}
where $\mc O^0_{r}$, $r\in\R[]$, represents a term that can be estimated by
\begin{equation}
\abs{\mc O^0_r}\le c\msp e^{rt}\qq\A\, 0\le t<\un.
\end{equation}
In inequality \re{6.80} we may choose $\de>0$ independently of $m\ge 1$. The terms inside the braces of that inequality are also an $\mc O^0_{-2\de}$ for an appropriate $\de>0$, and, because of \re{6.78},
\begin{equation}
\mc O^0_{-2\de}w=\mc O^0_{-\de}.
\end{equation}

Hence, applying the maximum principle to the function
\begin{equation}
w+\mu e^{-\de t}
\end{equation}
we derive an a priori estimate for $w$ by choosing $\mu$ large enough.

\cvm
(ii) It remains to  prove that the pointwise limit
\begin{equation}
\lim_{t\ra\un}(u(t,\xi)-2)e^{\frac tn}
\end{equation}
exists for any $\xi\in \Ss[n]$.

Using the scalar flow equation \fre{5.56} we deduce
\begin{equation}
\begin{aligned}
\dot{\tilde u}=\frac vFe^{\frac tn}+\tfrac1n\tilde u,
\end{aligned}
\end{equation}
where
\begin{equation}
\tilde u =(u-2)e^\frac tn
\end{equation}
and $F$ depends on
\begin{equation}
\check h^i_j=h^i_j+v^{-1}\tilde\vt \de^i_j.
\end{equation}

In view of the homogeneity of $F$ we further  conclude
\begin{equation}
\begin{aligned}
\dot{\tilde u}&=\frac v{F(\check h^i_je^{-\frac tn})}+\tfrac1n\tilde u\\
&=(-\tilde u)\Big\{\frac v{F(h^i_je^{-\frac tn}(-\tilde u)+v^{-1}\frac u{(1+\frac12 u)}\de^i_j)}-\tfrac1n\Big\}\\
&\ge -c \msp e^{-\frac tn}
\end{aligned}
\end{equation}
in view of our previous estimates, and we finally obtain
\begin{equation}
(\tilde u-nce^{-\frac tn})'\ge 0,
\end{equation}
from which the convergence result immediately follows. 
\ep


\providecommand{\bysame}{\leavevmode\hbox to3em{\hrulefill}\thinspace}
\providecommand{\href}[2]{#2}



\end{document}